\documentclass[12pt]{amsart}
\usepackage{amscd,amsmath,bm,latexsym,amssymb,amscd}
\usepackage[mathscr]{euscript} 
\usepackage[all]{xy}
\usepackage{layout}
\usepackage{pstricks}
\usepackage{enumerate}
\usepackage{dsfont}

\newtheorem{theorem}{Theorem}[section]
\newtheorem{corollary}[theorem]{Corollary}

\newtheorem{proposition}[theorem]{Proposition}
\theoremstyle{definition}
\newtheorem{definition}[theorem]{Definition}
\theoremstyle{definition}
\newtheorem{example}[theorem]{Example}
\theoremstyle{definition}

\newtheorem{remark}[theorem]{Remark}

\newtheorem*{theorem1*}{Theorem 3.1}
\newtheorem*{theorem2*}{Theorem 6.3}
\newtheorem*{theorem3*}{Theorem 2.2}
\newtheorem*{theorem4*}{Theorem 4.1}
\newtheorem*{theorem5*}{Theorem 7.2}
\newtheorem*{corollary6*}{Corollary 7.4}
\newtheorem*{corollary7*}{Corollary 7.5}
\newtheorem*{theorem8*}{Theorem 1.8}
\newtheorem*{theorem9*}{Theorem 6.5}
\def\Proof{\medskip\noindent{\bf Proof: }}

\def\Z{\mathbb{Z}}
\def\B{\mathbb{B}}
\def\En{\mathbb{E}}
\def\N{\mathbb{N}}
\def\C{\mathbb{C}}
\def\Q{\mathbb{Q}}
\def\R{\mathbb{R}}

\def\SS{\mathbb{S}}
\def\BC{B_{com}}
\def\EC{E_{com}}
\newcommand{\BONE}{\mathds 1}

\def\Bo{\mathcal{B}}
\def\Co{\mathcal{C}}
\def\Do{\mathcal{D}}
\def\E{\mathcal{E}}
\def\Fo{\mathcal{F}}
\def\H{\mathcal{H}}
\def\Go{\mathcal{G}}

\def\Ib{\mathbb{I}}
\def\Gr{\mathcal{G}r}

\def\Ma{\ \mathfrak{M}}
\def\O{\mathcal{O}}
\def\U{\mathcal{U}}

\def\Ro{\mathcal{R}}
\def\So{\mathcal{S}}
\def\Si{\textbf{S}}
\def\T{\mathcal{T}}

\def\i{\textbf{i}}
\def\j{\textbf{j}}
\def\g{\mathfrak{g}}

\def\t{\mathfrak{t}}

\def\z{\mathfrak{z}}
\def\x{{\bf{x}}}
\def\y{{\bf{y}}}

\def\colim{\displaystyle\mathop{\textup{colim}}}
\def\hocolim{\displaystyle\mathop{\textup{hocolim}}}
\DeclareMathOperator{\Hom}{\textup{Hom}}
\DeclareMathOperator{\maj}{\textup{maj}}
\DeclareMathOperator{\fmaj}{\textup{fmaj}}

\def\Tor{\textup{Tor}}
\def\Top{\textup{Top}}
\def\Ob{\textup{Ob}}
\def\Mor{\textup{Mor}}

\begin{document}

\title[A Classifying Space for Commutativity in Lie Groups] 
{A Classifying Space for Commutativity in Lie Groups}

\author[A.~Adem]{Alejandro Adem}
\address{Department of Mathematics,
University of British Columbia, Vancouver BC V6T 1Z2, Canada}
\email{adem@math.ubc.ca}

\author[J.~M.~G\'omez]{Jos\'e Manuel G\'omez}
\address{Department of Mathematics,
Johns Hopkins University, Baltimore, MD 21218, USA.\newline  
Departamento de Matem\'aticas,
Universidad Nacional de Colombia, Medellin, AA 3840, Colombia}
\email{jmgomez0@unal.edu.co}
\thanks{The first author was supported by NSERC.
The second author would like to thank PIMS for hosting him 
when part of this work 
was completed.}

\begin{abstract}
In this article we consider a space $\BC G$ assembled from
commuting elements in a Lie group $G$ first defined in \cite{ACT}. 
We describe homotopy-theoretic
properties of these spaces using homotopy colimits,
and their role as a classifying space for 
\textsl{transitionally commutative} bundles. We prove that 
$\Z\times \BC U$
is a loop space and define a notion of commutative K--theory
for bundles over a finite complex $X$ which is isomorphic to
$[X,\Z\times\BC U]$.
We compute the rational
cohomology of $\BC G$ for $G$ equal to any of the classical
groups  $SU(r)$, $U(q)$ and $Sp(k)$, and exhibit the 
rational cohomologies of $\BC U$, $\BC SU$ and $\BC Sp$ as
explicit polynomial rings.
\end{abstract}

\maketitle

\section{Introduction}
Let $G$ denote a topological group and consider the spaces 
$\{\Hom(\Z^n,G)\}_{n\ge 0}$ of ordered
commuting $n$-tuples in $G$. In \cite{ACT} it was shown that they
can be assembled into a simplicial space 
where the resulting geometric realization,
denoted here by $\BC G$, is the first term in an increasing filtration
of $BG$.
The universal bundle over $BG$ pulls  
back to a principal bundle over $\BC G$ with total space
$\EC G$ that can 
also be described  simplicially 
(see \S \ref{spaces B_{com}G}). 
In this paper we analyze the properties of $\EC G$ and $\BC G$ 
when $G$ is a Lie group. We also study
variants of our constructions denoted $\EC G_{\BONE}$ and 
$\BC G_{\BONE}$ which arise from the components of the identity 
in the commuting varieties (see \S \ref{properties} for details).
These two constructions agree if the spaces $\Hom(\Z^n,G)$ are path--connected
for every $n\ge 0$. It can be shown that for
a compact Lie group $G$, this condition
is equivalent to the property that  
the maximal abelian subgroups of $G$ are precisely the maximal tori
(see the proof of \cite[Proposition 2.5]{AG}).
For example, this condition holds for the classical groups 
$SU(r)$, $U(q)$ and $Sp(k)$ (see \cite[Theorem 5.2]{Borel}), 
and therefore
for any of their finite cartesian products.
 
We start by applying the recent work in \cite{PS} to reduce
matters to compact Lie groups:

\begin{theorem1*}
If $G$ is a real or complex reductive algebraic group
with maximal compact subgroup $K$, then the inclusion map
$K\subset G$ induces homotopy equivalences
$\BC K\simeq\BC G$ and $\EC K\simeq\EC G$.
\end{theorem1*}

A similar statement is true for the variants $\EC G_{\BONE}$ and 
$\BC G_{\BONE}$ (see \S \ref{properties} for details). 
Based on this we can focus on the case of a compact connected 
Lie group $G$.
The connected component of the identity 
in the commuting variety $\Hom(\Z^{n},G)$ has the key feature that 
any $n$-tuple in it can be conjugated into a maximal torus in $G$.
Using this we obtain a  natural identification  
$\BC G_{\BONE}\cong\colim_{S\in \T(G)}BS$, where $\T(G)$ is
the topological poset formed by the maximal tori and their
intersections under inclusion (see \ref{poset T(G)}).
We describe the homotopy of these spaces using  
more tractable homotopy colimits defined over a discrete 
category.
Let $Z=Z(G)$ be the center of $G$ and write 
$n={\rm{rank}}(G)-{\rm{rank}}(Z)\ge 0$. Consider the poset 
$\So(n)$  
consisting of all the nonempty subsets of $\{0,1,\dots, n\}$, with 
the order given by the reverse inclusion of sets.   
For each $G$ we have  functors $\Fo_{G},\H_{G}:\So(n)\to \Top$  
such that the following holds (see Section 6 for the definitions of 
$\Fo_{G}$ and $\H_{G}$).

\begin{theorem2*}\label{theorem 6.3}
Suppose that $G$ is a compact, connected Lie group. Then there
is a natural homotopy equivalence 
$\hocolim_{\i\in \So(n)}\Fo_{G}(\i)\simeq  B_{com}G_{\BONE}$.
\end{theorem2*}

\begin{theorem9*}\label{theorem 6.5}
Suppose that $G$ is a compact connected Lie group. Then there is a 
natural $G$-equivariant homotopy equivalence
$\hocolim_{\i\in \So(n)}\H_{G}(\i)\simeq  E_{com}G_{\BONE}$.
\end{theorem9*}

In terms of bundle theory, we prove that $\BC G$ is a classifying
space for bundles which are \textsl{transitionally commutative}.

\begin{theorem3*}\label{geometric interpretation}
Suppose that $G$ is a Lie group and let $f:X\to BG$ denote the 
classifying map of a  principal
$G$--bundle  $q:E\to X$  over the finite CW--complex $X$. 
Then up to homotopy, $f$ factors through $\BC G$ if and only if there 
is an open cover of $X$ on which the bundle is trivial over each open 
set and such that on intersections the transition functions commute 
when they are simultaneously defined.
\end{theorem3*} 

From this we can define the notion of equivalence between 
transitionally commutative vector bundles and thus define 
\textsl{commutative} K--theory $K_{com}(X)$ in a manner
analogous to ordinary complex K--theory. Let
$U=\colim_{n\to\infty}U(n)$, then we can
establish that $\BC U$ plays a role similar to $BU$:

\begin{theorem4*}
The space $\Z\times \BC U$ is a loop space and for any finite
$CW$--complex $X$ there is a natural
isomorphism of groups $K_{com}(X)\cong [X, \Z\times \BC U]$.
\end{theorem4*}

In \cite{AGLT} it is proved that  $\Z\times \BC U$ 
is in fact an infinite loop space so that commutative K-theory 
forms part of a generalized cohomology theory. 
Having established the role played by $\BC G$ in bundle theory,
it seems natural to compute its cohomology. As can be seen from
the homotopy colimit model and computations for $SU(2)$ 
(see Example \ref{case SU(2)}), we expect
these spaces to have rather intricate torsion.
Here we focus on calculations for the rational cohomology 
(inverting the order of the Weyl group $W$ would suffice).

\begin{theorem5*}\label{free module}
Suppose that  $G$ is a compact, connected Lie group.
Then 
$H^{*}(\BC G_{\BONE};\Q)$ is a free module over $H^{*}(BG;\Q)$ 
of rank $|W|$, where $W$ is the corresponding Weyl group.
\end{theorem5*}

On the other hand, by \cite[Theorem 6.1]{ACT} we have a natural 
isomorphism 
$H^*(\BC G_{\BONE};\Q)\cong (H^*(G/T;\Q)\otimes H^*(BT;\Q))^W$, where
the Weyl group $W$ acts diagonally on the tensor product. As a 
corollary we deduce the following algebra isomorphism:

\begin{corollary6*}
Suppose that $G$ is a connected compact Lie group with
maximal torus $T$ and associated Weyl group $W$.
Then there is a natural 
isomorphism of rings 
\[
H^{*}(E_{com}G_{\BONE})
\cong (H^*(G/T)\otimes H^*(G/T))^W,
\]
and the Poincar\'e series of $\BC G_{\BONE}$ and 
$E_{com}G_{\BONE}$ satisfy 
$P_{\BC G_{\BONE}}(t)=P_{BG}(t)P_{E_{com}G_{\BONE}}(t)$.
\end{corollary6*}

From this the following is derived.

\begin{corollary7*} The following three statements are equivalent for a compact
connected Lie group $G$:
(1)~$E_{com}G_{\BONE}$ is contractible;
(2)~$E_{com}G_{\BONE}$ is rationally acyclic; and
(3)~$G$ is abelian.
\end{corollary7*}

Using the theory of multisymmetric polynomials in 
\S 8 we 
provide combinatorial descriptions and Poincar\'e series for these 
algebras in the case of the
classical groups $SU(r)$, $U(q)$ and $Sp(k)$. Taking
limits we obtain that the algebras $H^*(\BC U;\Q)$, $H^*(\BC SU;\Q)$ and 
$H^*(\BC Sp;\Q)$ are polynomial algebras on countably 
many generators (Corollaries 8.3, 8.5 and 8.9 respectively).
For example
we have an isomorphism of $\Q$-algebras
 \[
 H^{*}(B_{com}U;\Q)\cong  \Q[z_{a,b}~|~ (a,b)\in\N^{2} \text{ and } b>0],
 \]
where the elements $z_{a,b}$ are polynomial generators
of degree $2a + 2b$.

\medskip

This paper is organized as follows: in \S 2 we describe basic
properties of the spaces $\EC G$ and $\BC G$ for $G$ a topological
group; in \S 3 we focus
on the case when $G$ is a Lie group; in \S 4 commutative K-theory 
is introduced; \S 5 describes the topological poset generated by the 
maximal tori in a Lie group $\T(G)$; 
in section \S 6  we derive the decompositions of 
$B_{com}G_{\BONE}$ and $E_{com}G_{\BONE}$ as 
homotopy colimits; \S 7 deals with cohomology calculations; 
in \S 8 we consider the particular cases when $G=SU(n)$, 
$U(n)$ and $Sp(n)$ and finally \S 9 is an  appendix 
where it is proved that $[B_{com}G]_{*}$ is a proper simplicial space 
for any Lie group $G$. 
We are grateful to the referee for providing helpful comments.

\section{Definitions and basic properties of the spaces 
$B_{com}G$  and $E_{com}G$}\label{spaces B_{com}G}

In this section we study general properties of the spaces $B_{com}G$ 
and $E_{com}G$ which are constructed by assembling the different 
spaces of ordered commuting $k$-tuples in a topological group 
$G$. These spaces were first introduced in \cite{ACT} where
their basic properties were derived, mostly for the case of finite groups. 

\medskip

Suppose that $G$ is a topological 
group. For technical reasons we will assume that $G$ is locally
compact, Hausdorff and that $1_{G}\in G$ is a non-degenerate 
basepoint.  We can associate to $G$ a simplicial space, 
denoted by $[B_{com}G]_{*}$, in the following 
way. For any integer $n\ge 0$ define 
\[
[B_{com}G]_n:=\Hom(\Z^{n},G)\subset G^n.
\] 
Note that $[B_{com}G]_n$ can be identified with the subset of $G^{n}$ 
consisting of all ordered 
commuting $n$-tuples in the group $G$, and as such it is given 
the subspace topology. The face and degeneracy maps are 
defined by
\[
s_{j}(g_{1},...,g_{n})=(g_{1},...,g_{j},1_{G},g_{j+1},...,g_{n}),
\] 
and
\begin{equation*}
\partial _{i}(g_{1},...,g_{n})=\left\{
\begin{array}{ccc}
(g_{2},...,g_{n}) & \text{if } &i=0, \\
(g_{1},...,g_{i}g_{i+1},...,g_{n}) & \text{if } &0<i<n, \\
(g_{1},...,,g_{n-1}) & \text{if } &i=n.
\end{array}
\right.
\end{equation*}
The different $s_{i}$'s and $\partial_{j}$'s are 
well defined and satisfy the simplicial identities as these maps 
are precisely the restrictions of the degeneracy and face maps in 
the bar construction $[BG]_{*}$. We denote by $B_{com}G$ the 
geometric realization of the simplicial space $[B_{com}G]_*$. 
As shown in \cite{ACT}, the space $B_{com}G$ 
is in fact the first space in an increasing filtration of the classifying 
space of $G$ defined using the descending central series of the free
groups.   Similarly we can define 
$[E_{com}G]_n:=\Hom(\Z^{n},G)\times G\subset G^{n+1}$
and use the analogous face and degeneracy maps to define a 
simplicial space $[E_{com}G]_*$ and its geometric 
realization $E_{com}G$.

The projection on the first $n$-coordinates 
$[E_{com}G]_*\to [B_{com}G]_*$ defines a simplicial map and 
therefore at the level of geometric realizations  
we obtain a continuous 
map $p_{com}:E_{com}G\to B_{com}G$. This defines a principal 
$G$-bundle that can be seen as the restriction of the universal 
principal $G$-bundle $p:EG\to BG$, and we have  
a morphism of principal $G$-bundles 
\begin{equation*}
\begin{CD}
E_{com}G@> >> EG\\
@V{p_{com}}VV     @VV{p}V\\
B_{com}G@>{i}>>
BG.\\
\end{CD}
\end{equation*}
Note that up to homotopy this gives rise to
a fibration sequence 
$E_{com}G\to B_{com}G\to BG$.

Recall that the 
bundle $p:EG\to BG$ is universal in the sense 
that if $q:E\to X$ is a principal $G$-bundle 
over a CW--complex $X$, then we can find a 
continuous map $f:X\to BG$ such 
that $q:E\to X$ is isomorphic to $f^{*}{p}:f^{*}(EG)\to X$. 

\begin{definition}
Suppose that $X$ is a CW--complex. We say that a 
principal $G$-bundle $q:E\to X$ is \textsl{transitionally commutative} 
if and only if we can find an open cover $\{U_{i}\}_{i\in I}$ of $X$ 
such that the bundle $q:E\to X$ is trivial over each $U_{i}$ 
and the transition functions $\rho_{i,j}:U_{i}\cap U_{j}\to G$ commute 
with each other whenever they are simultaneously defined. 
\end{definition}

When $G$ is a Lie group the bundle $p_{com}:E_{com}G\to B_{com}G$ is 
a universal bundle for transitionally commutative principal $G$-bundles.
To make this precise we need to establish the following notation.
Let $k\ge 0$ be an integer and consider the standard $k$-simplex 
\[
\Delta_{k}=\{(t_{0},\dots, t_{k})\in \R^{k+1} ~|~ t_{i}\ge 0, 
\sum_{j=0}^{k}t_{j}=1 \}.
\]
Note that the symmetric group $\Sigma_{k+1}$ acts by permutation 
on the vertices of $\Delta_{k}$ and this action can be extended 
to a linear action on $\Delta_{k}$.
Suppose that we have a sequence of integers $\i:=\{0\le i_{1}<\cdots <i_{q}\le k\}$ 
and let $e_{\i}$ denote the element in 
$\Delta_{k}$ given by  
$e_{\i}=\delta_{i_{q}}\cdots \delta_{i_{1}}(\frac{1}{k-q+1},\dots,\frac{1}{k-q+1})$, 
where $\delta_{i_{1}},\dots,\delta_{i_{q}}$ denote the different face maps. 
In other words, the element $e_{\i}\in \Delta_{k}$ has barycentric 
coordinates $(t_{0},\dots, t_{k})$ given by 
\begin{equation*}
t _{j}=\left\{
\begin{array}{ccc}
0 & \text{if } &j\in \{i_{1},\dots,i_{q}\}, \\
\frac{1}{k-q+1} & \text{if } &j\notin \{i_{1},\dots,i_{q}\}.
\end{array}
\right.
\end{equation*} 
The length of such a sequence $\i$ is defined to be the number  
$|\i|=q$.
Observe that if $\i$ runs through all the different sequences 
of integers $\i:=\{0\le i_{1}<\cdots <i_{q}\le k\}$ then the collection 
$\{e_{\i}\}_{\i}$ is precisely the collection of barycenters associated 
to all the non-empty faces in the standard simplex in $\Delta_{k}$.   
Let $\{W^{k}_{q} \}_{q=0}^{k}$ be a collection 
of open sets in $\Delta_{k}$ satisfying the following 
properties:
\begin{enumerate}\label{simplex cover}
\item for each $0\le q\le k$ and each $\i=\{0\le i_{1},\dots, i_{q}\le k \}$ 
there is an open neighborhood  $V_{\i}^{k}$ of $e_{\i}$ in such a way that  
$W^{k}_{q}= \bigsqcup_{|\i|=q+1}V_{\i}^{k}$,

\item if $t=(t_{0}, \dots, t_{k})\in V_{\i}^{k}$ for some 
$\i=\{0\le i_{1}<\cdots <i_{q}\le k\}$, then $t_{j}>0$ 
if $j\notin \{i_{1},\dots, i_{q}\}$,

\item each open set $W^{k}_{q}$ is invariant under the action of $\Sigma_{k+1}$, and 

\item the sets $\{W^{k}_{q} \}_{q=0}^{k}$ form an open cover of $\Delta_{k}$.
\end{enumerate}
Clearly such an open cover exists and can be constructed in an inductive way.
Using this we have the following geometric 
description for the bundle 
$p_{com}:E_{com}G\to B_{com}G$.

\begin{theorem}
Suppose that $G$ is a Lie group and let $f:X\to BG$ denote the 
classifying of a  principal
$G$--bundle  $q:E\to X$  over the finite CW--complex $X$. Then up to homotopy,
$f$ factors through $\BC G$ if and only 
if $q$ is transitionally commutative.
\end{theorem}

\Proof Let  $X$ be a finite CW--complex.  Assume that the 
classifying map of $q:E\to X$ factors through $B_{com}G$; that is, 
suppose that the classifying map of this bundle is of the form 
$f:X\to B_{com}G\subset BG$. We will show first 
that $q$ is transitionally commutative. 
By Proposition \ref{proper simplicial} in the appendix we have that  
$[B_{com}G]_{*}$ is a proper simplicial space (see the appendix for the 
definition of a proper simplicial space). It follows that 
the geometric realization of $[B_{com}G]_{*}$ 
is equivalent to Segal's fat geometric 
realization where the equivalences associated to the 
degeneracy maps are ignored (see \cite[Appendix A]{Segal}). 
This realization is denoted here  by $\B_{com}G$, and similarly we
have $\En_{com}G$.
For the first part of the proof it will be more convenient for us to
work with the (equivalent) principal $G$--bundle 
$p_{com}:\En_{com}G\to \B_{com}G$.

By definition
\[
\B_{com}G:=\left(\bigsqcup_{n\ge 0}\Hom(\Z^{n},G)\times \Delta_{n}\right)/\sim,
\]
where $((g_{1},\dots,g_{k}),\delta_{i}u)\sim 
(\partial_{i}(g_{1},\dots,g_{k}),u)$ with $u\in \Delta_{k-1}$.  
For each $k\ge 0$, let 
\[
F_{k}\B_{com}G=\text{Im}\{\bigsqcup_{0\le n\le k}
\Hom(\Z^{n},G)\times \Delta_{n}\} \subset \B_{com}G. 
\]
In this way we obtain an increasing filtration of $\B_{com}G$
\[
F_{0}\B_{com}G\subset F_{1}\B_{com}G\subset\cdots\subset 
F_{k}\B_{com}G\subset \cdots \subset \B_{com}G   
\]
and $\B_{com}G=\colim_{k\to \infty}F_{k}\B_{com} G$. Since $X$ is a finite 
CW--complex, then we can find some $k\ge 0$ such that the map 
$f$ factors through $F_{k}\B_{com}G$; that is,  
$f:X\to F_{k}\B_{com} G\subset \B_{com}G$. 
It suffices to show the result for the restriction of the 
universal principal $G$-bundle $p_{com}:\En_{com} G\to \B_{com}G$ over $F_{k}\B_{com}G$ 
for each $k\ge 0$ fixed.  For this fix $\i:=\{0\le i_{1}<\cdots <i_{q}\le k\}$
a sequence of integers and assume that 
$(g_{1},\dots,g_{k})\in \Hom(\Z^{k},G)=[B_{com}G]_{k}$.  Then we can see 
$(g_{1},\dots,g_{k},1)\in \Hom(\Z^{k},G)\times G=[E_{com}G]_{k}$ and we 
define 
\[
\varphi_{\i}(g_{1},\dots, g_{k}):=
\pi_{k-q}(\partial_{i_{1}}\cdots \partial_{i_{q}}(g_{1},\dots,g_{k},1))^{-1}.
\] 
In the above equation $\partial_{i_{1}},\dots,\partial_{i_{q}}$ denote 
the face maps in the simplicial space $[E_{com}G]_{*}$ and 
$\pi_{n-k}:\Hom(\Z^{n-k-1})\times G\to G$ denotes the projection onto 
the last coordinate. For example if $(g_{1},g_{2}, g_{3})\in \Hom(\Z^{3},G)$ 
and $\i=\{2,3\}$, then 
$\varphi_{\i}(g_{1},g_{2}, g_{3})=g_{3}^{-1}g_{2}^{-1}$ and if 
$\j=\{1,3\}$ then $\varphi_{\i}(g_{1},g_{2}, g_{3})=g_{3}^{-1}$.
 Note that in general 
for any $(g_{1},\dots, g_{k})\in \Hom(\Z^{k},G)$ and any sequence $\i$ we have 
$\varphi_{\i}(g_{1},\dots, g_{k})=g_{j_{1}}^{-1}\cdots g_{j_{r}}^{-1}$ for 
suitable integers $j_{1},\dots, j_{r}$. The different functions $\varphi_{\i}$
can be used to define local sections of the restriction of the bundle 
$p_{com}:\En_{com}G\to \B_{com}G$ over $F_{k}\B_{com}G$. Indeed, for 
each $0\le q\le k$ let $U_{q}^{k}$ be the image of  
$\Hom(\Z^{k},G)\times W^{k}_{q}$ in $F_{k}\B_{com}G$. 
Defined this way, each $U_{q}^{k}$ is an open set in $F_{k}\B_{com}G$ and 
the collection $\{U_{q}^{k}\}_{q=0}^{k}$ forms an open cover of 
$F_{k}\B_{com} G$. Define $\sigma_{q}:U_{q}^{k}\to p_{com}^{-1}(U_{q}^{k})$ 
in the following way. Let $x\in U_{q}^{k}$ and 
write $x=[(g_{1},\dots,g_{k}),t]$ for some 
$(g_{1},\dots, g_{k})\in \Delta_{k}$ and $t\in W_{q}^{k}$. We define 
$\sigma_{q}(x):=[(g_{1},\dots,g_{k},\varphi_{\i}(g_{1},\dots,g_{k})),t]$, 
provided that $t\in V_{\i}^{k}$. The functions $\varphi_{\i}$ are 
defined precisely so that the function $\sigma_{q}$ is well defined 
and continuous over $U_{q}^{k}$. Therefore  $\sigma_{q}$ is  
a continuous section of the restriction of $p_{com}$ 
over $U^{k}_{q}$ making it a trivial principal $G$-bundle. 
Moreover, with the trivializations provided by these 
sections, if $x=[(g_{1},\dots,g_{k}),t]\in U^{k}_{r}\cap U^{k}_{q}$, 
then transition function 
$\rho_{r,q}:U^{k}_{r}\cap U^{k}_{q}\to G$ is such that 
\[
x=[(g_{1},\dots,g_{k}),t]\mapsto g_{j_{0}}^{\pm 1}\cdots g_{j_{r}}^{\pm 1},
\]
for a suitable sequence of integers $j_{0},\dots, j_{r}$.
In particular it follows that the different 
transition functions $\rho_{r,q}$ are pairwise commutative 
whenever they are simultaneously defined as 
$(g_{1},\dots,g_{k})\in \Hom(\Z^{k},G)$ for any 
$x=[(g_{1},\dots,g_{k}),t]\in U^{k}_{r}\cap U^{k}_{s}$.

Conversely, suppose that $q:E\to X$ is a 
transitionally commutative principal $G$-bundle.
Then we can find an open cover  
$\U:=\{U_{i}\}_{i\in I}$ of $X$ such that 
the bundle $q:E\to X$ is trivial over each $U_{i}$ 
and the transition functions $\rho_{i,j}:U_{i}\cap U_{j}\to G$ commute 
with each other whenever they are simultaneously defined. By passing 
to a refinement of $\U$, we can assume that each nonempty 
intersection of the sets $U_{i}$'s is contractible. Moreover, this 
cover can be reduced to a countable cover. Let 
$\U=\{U_{i}\}_{i\ge 0}$ be the resulting open cover of $X$. 
For each $i\ge 0$, fix $\varphi_{i}:q^{-1}(U_{i})\to U_{i}\times G$ a 
trivialization of the restriction of $q$ over $U_{i}$. These 
trivializations  define transition 
functions $\rho_{i,j}:U_{i}\cap U_{j}\to G$  whenever 
$U_{i}\cap U_{j}\ne \emptyset$ and satisfy the 
cocycle condition  $\rho_{ik}=\rho_{ij} \rho_{jk}$ whenever they are 
defined and are pairwise commutative by assumption. 
Consider the nerve $N_{*}(\U)$ of the cover $\U$. This is 
a simplicial set with 
$N_{k}(\U)=\bigsqcup_{0\le i_{1}\le \cdots\le i_{k+1}} 
U_{i_{1}}\cap\cdots\cap U_{i_{k+1}}$.
The different transition functions can be used to define a map of 
simplicial spaces 
$\rho_{*}:N_{*}(\U)\to [B_{com}G]_*$
in the following way. Suppose that 
$x\in U_{i_{1}}\cap\cdots\cap U_{i_{k+1}}$ for some 
$0\le i_{1}\le \cdots\le i_{k+1}$. Define 
$\rho_{k}(x)=(\rho_{i_{1}i_{2}}(x),\rho_{i_{2}i_{3}}(x),\dots,
\rho_{i_{k}i_{k+1}}(x))\in \Hom(\Z^{k},G)$.
It is easy to see that this defines a map of simplicial spaces and in 
particular it induces a continuous map $g=|\rho_{*}|:N(\U)\to \BC G$. 
Since the cover $\U$ was chosen so that 
each nonempty intersection of sets in 
$\U$ is contractible, then the natural map $\alpha:N(\U)\to X$ is a 
homotopy equivalence (see for example \cite[Corollary 4G.3]{Hatcher}). 
Let $\beta:X\to N(U)$ be a homotopy inverse 
of $\alpha$. Then $f:=g\circ\beta:X\to \BC G$ is a continuous 
map that classifies the principal $G$-bundle  $q:E\to X$.
\qed

\medskip

As a consequence of the proof of the previous theorem, we have that 
the restriction of the bundle $p_{com}:\En_{com}G\to \B_{com}G$ 
to each $F_k\B_{com}G$ defines a transitionally commutative 
principal $G$-bundle.  From this we infer that the 
bundle $p_{com}:\En_{com}G\to \B_{com}G$ is itself transitionally commutative,
as is the equivalent bundle $p_{com}:E_{com}G\to B_{com}G$. 

As an application of the previous theorem suppose that $G$ 
is a Lie group and that $X$ is a finite CW--complex 
for which we can find an open cover $X=U\cup V$ 
with both $U$ and $V$ contractible.  Let $q:E\to X$ 
be any principal $G$-bundle over $X$. Then the restriction of $q$ 
over $U$ and $V$ is trivial since $U$ and $V$ are contractible.
Over this trivialization there is only one transition function and thus 
any such principal $G$-bundle over $X$ is transitionally commutative. 
By the previous theorem we conclude that the classifying 
map of the bundle $q:E\to X$ factors through $\BC G$ up to homotopy.  
This situation applies in particular to $X=\SS^{n}$ for any $n\ge 0$. 
Therefore the inclusion map $i:\BC G\hookrightarrow  BG$ induces 
a surjective map
$i_{\#}:[\SS^{n},\BC G]\to [\SS^{n},BG]$. The following corollary
is an immediate consequence after modifying for basepoints and 
using the fibration $E_{com}G\to B_{com}G\to BG$.

\begin{corollary}\label{ses homotopy} Let $G$ be a Lie group, 
then the map $i:\BC G\to BG$ induces a 
surjection $i_{*}:\pi_{n}(\BC G)\to \pi_{n}(BG)$ for every $n\ge 0$ 
and in particular, for every $n\ge 0$ we have a short exact sequence 
$1\to \pi_{n}(E_{com}G)\to \pi_{n}(\BC G)\stackrel{i_{*}}
{\rightarrow} \pi_{n}(BG)\to 1$.
\end{corollary}

\begin{remark}
Suppose that $G$ is a connected Lie group. Then by 
\cite[Theorem 6.3]{ACT}  the fibration sequence  
$\Omega E_{com}G\to \Omega B_{com}G\to \Omega BG $
has a natural continuous section 
$\sigma(G):\Omega BG \to  \Omega B_{com}G$.  
This implies that for such groups and every $n\ge 0$ 
the short exact sequence in homotopy groups 
obtained in the previous corollary splits naturally. 
Moreover, by \cite[Theorem 6.3]{ACT}  
there is a natural homotopy equivalence
$\theta(G):G\times \Omega E_{com}G\to \Omega B_{com}G$.
\end{remark}

\medskip

Suppose that $X$ is a finite CW--complex. Note that 
two principal $G$-bundles  $q_{0}:E_{0}\to X$ 
and  $q_{1}:E_{1}\to X$ are isomorphic if and only if we can 
find a  principal $G$-bundle $p:E\to X\times [0,1]$ such that 
$q_{0}=p_{|p^{-1}(X\times\{0\})}$ and $q_{1}=p_{|p^{-1}(X\times\{1\})}$. 
Suppose now that $q_{0}:E_{0}\to X$ and  $q_{1}:E_{1}\to X$ are 
two transitionally commutative principal $G$-bundles. Then we 
say that these bundles are \textit{transitionally commutative 
isomorphic} if we can find a  transitionally 
commutative principal $G$-bundle $p:E\to X\times [0,1]$ such that 
$q_{0}=p_{|p^{-1}(X\times\{0\})}$ and $q_{1}=p_{|p^{-1}(X\times\{1\})}$. 
Thus we can identify the set $[X,\BC G]$ with the set of transitionally 
commutative isomorphism classes of transitionally commutative 
principal $G$-bundles over $X$. In other words, the space 
$B_{com}G$ is a classifying space for transitionally commutative 
bundles. If two transitionally 
commutative principal $G$-bundles are transitionally 
commutative isomorphic then they are isomorphic as principal 
$G$-bundles. However, the converse is not true as is 
demonstrated in the next example. 

\begin{example}\label{nontrivial class SU(2)}
Let $G=SU(2)$ and $T\subset G$ the maximal torus, which in
this case is a circle. The quotient $G/T$ can be identified
with the sphere $\SS^2$, let $f:\SS^2\to G/T$ be a fixed
homeomorphism. Now the action map $G\times T^n \to \Hom(\Z^n,G)$ 
defined by $(g,t_1,\dots ,t_n)\mapsto (gt_1g^{-1},\dots ,gt_ng^{-1})$
factors through $G/T\times T^n$. Looking at the realizations
of the respective simplicial spaces, this defines a
map $\theta: G/T\times BT\to B_{com}G$. According to \cite{ACT}, 
Theorem 6.1, this gives rise to a rational cohomology isomorphism 
$H^*(B_{com}G,\Q)\to H^*(G/T\times BT,\Q)^{W}$, where
$W= \Z/2\Z$ is the Weyl group. This group acts
through the sign representation both on the generator 
$a\in H^2(BT,\Q)$ and on the top class $b\in H^2(G/T,\Q)$.
The invariant classes $a^2$ and $ab$ correspond to a basis 
for $H^4(B_{com}G,\Q)$. Now let 
$g:\SS^2\to BT\cong \C P^{\infty}$ be a representative of 
a generator of 
$\pi_{2}(\C P^{\infty})\cong \Z$. Consider the map 
$h:\SS^{2}\times \SS^{2}\to B_{com}SU(2)$ given by the
composition $\theta\circ (f\times g)$. In the next section
we shall see that $B_{com}G$ is 3--connected for $G=SU(2)$,
which implies that
this map is null homotopic on $\SS^2\vee\SS^2$ and so
defines a map $\tilde{h}:\SS^4\to B_{com}G$. By construction
this map is non-trivial in rational cohomology, corresponding
to the element $ab$. 
Moreover, if $i:B_{com}G\to BG$ denotes the inclusion 
map, then the composition $i\circ \tilde{h}$ is trivial in cohomology 
(as the Chern class in dimension four corresponds to $b^2$) and so 
is null homotopic. It follows that the principal $G$-bundle
over $\SS^{4}$ induced by $\tilde{h}$ is trivial as a principal
$G$-bundle but not as a transitionally
commutative principal $G$-bundle.
\end{example}

\begin{remark}
If $G=U(k)$, with $k>1$, then the map 
$i_{*}:\pi_{n}(B_{com}U(k))\to \pi_{n}(BU(k))$ 
cannot be an isomorphism for every $n\ge 0$. If this were true 
$i$ would be a homotopy equivalence; however it follows from
\cite{ACT}, Theorem 6.1 that the rational cohomology of $B_{com}U(k)$ 
is not isomorphic to that of $BU(k)$. 
More generally we shall see that if $G$ is a compact connected Lie
group which is not abelian, then $E_{com}G$ 
is \textbf{not contractible}, unlike the classical universal space $EG$.
\end{remark}

\section{Properties of $\BC G$ and $E_{com}G$ when 
$G$ is a Lie Group}\label{properties}

In this section we focus our attention on the case when 
$G$ is a real or complex reductive algebraic group. 
We can consider $G$ as a real or complex Lie group respectively. 
Let $K\subset G$  be a maximal compact subgroup; it is well 
known that such a group always exists and 
the inclusion map $i:K\hookrightarrow G$ is a strong deformation 
retract. However, in general there is no retraction $r:G\to K$ 
that preserves commutativity; see for example \cite{Souto} where 
the nonexistence of such a retraction was  
proved for the groups $SL_{n}(\C)$ with $n\ge 8$. On the 
other hand, by  \cite[Corollary 1.2]{PS} the
inclusion  $\Hom(\Z^{n},K)\hookrightarrow \Hom(\Z^{n},G)$ 
is a strong deformation retract. We show here that this can be used 
to prove that the inclusion $i:B_{com}K\hookrightarrow \BC G$ is also 
a strong deformation retract. 

\begin{theorem}\label{reduction to compact}
Suppose that $G$ is a real or complex reductive algebraic group
and let $K$ be a maximal compact subgroup. Then 
the inclusion map $i:K\hookrightarrow G$ induces homotopy 
equivalences $i:B_{com}K\to \BC G$ and $i:E_{com}K\to E_{com}G$.
\end{theorem}

\Proof
We can see $G$ as a (real or complex) 
Lie group and thus  by Proposition \ref{proper simplicial} in the 
appendix $[B_{com}G]_{*}$ 
and $[B_{com}K]_{*}$ are proper simplicial spaces. 
The inclusion map $i:K\hookrightarrow G$ induces a map of 
simplicial spaces  $i_{*}:[B_{com}K]_*\to [B_{com}G]_*$ that is a 
level-wise homotopy equivalence. 
By \cite[Theorem A.4]{May1} we conclude that the induced map at the 
level of geometric realizations $i:B_{com}K\to \BC G$ is a  homotopy 
equivalence. Next we prove that $i:E_{com}K\to E_{com}G$ 
is a homotopy equivalence. For this, note that the inclusion map 
$B_{com}K\to \BC G$ induces a morphism 
of the corresponding fibrations 
\begin{equation*}
\begin{CD}
E_{com}K@> >> B_{com}K@>>> BK\\
@V{i}VV  @V{i}VV     @V{i}VV\\
E_{com}G@> >> \BC G@>>> BG.\\
\end{CD}
\end{equation*} 
This diagram induces a commutative diagram between the 
corresponding long exact sequences in homotopy groups. 
Since the inclusion maps 
$i:B_{com}K\to \BC G$ and 
$i:BK\to BG$ are homotopy equivalences, 
by the five lemma it follows that the inclusion map 
$i:E_{com}K\to  E_{com}G$
is also a homotopy equivalence.
\qed

\medskip

If $G$ is a compact Lie group then it can be given 
the structure of a real algebraic variety that is reductive by 
complete reducibility. Let $G_{\C}$ 
denote its complexification. Then by the previous theorem it follows 
that $\BC G$ and $B_{com}G_{\C}$ are homotopy equivalent and 
similarly for $E_{com}G$ and $E_{com}G_{\C}$.  This 
shows that we can work in the category of compact Lie groups 
without loss of generality whenever we want to study the spaces 
$\BC G$ and $E_{com}G$ for a real or complex reductive 
algebraic group $G$.

Suppose that $G$ is a topological group; 
$\Hom(\Z^{n},G)$ may fail to be path--connected even if we assume 
that $G$ is path--connected or simply connected. 
For every $n\ge 0$ define
$\Hom(\Z^{n},G)_{\BONE}$ to be the path--connected component of 
$\Hom(\Z^{n},G)$ containing the trivial representation 
$\BONE:\Z^{n}\to G$.  It is easy to see that the collection 
$\{\Hom(\Z^{n},G)_{\BONE}\}_{n\ge 0}$ forms a simplicial subspace  of
$[B_{com}G]_*$. We denote by $\BC G_{\BONE}$ its geometric realization.
When $G$ is a compact Lie group the path--connected component 
$\Hom(\Z^{n},G)_{\BONE}$ has the following important 
feature as already pointed out in \cite[Lemma 4.2]{Baird}. An $n$-tuple 
$(g_{1},\dots,g_{n})$ of elements in $G$ belongs to 
$\Hom(\Z^{n},G)_{\BONE}$ if and only if there is a maximal 
torus $T\subset G$ that contains $g_{1},\dots,g_{n}$. On the other 
hand, if $G$ is a complex reductive algebraic variety then a 
commuting tuple $(g_{1},\dots, g_{n})$ belongs to 
$\Hom(\Z^{n},G)_{\BONE}$ if and only if there is a torus $T\subset G$
containing the semisimple part of the Jordan decomposition of $g_{i}$ 
for all $1\le i\le n$.  The spaces $\BC G$ and $\BC G_{\BONE}$ agree if 
$\Hom(\Z^{n},G)$ is path--connected for all $n\ge 0$. When $G$ is a 
compact Lie group this  
is the case if and only a subgroup $A\subset G$ is a
maximal abelian subgroup in $G$ if and only if $A$
is a maximal torus in $G$ by \cite[Proposition 2.5]{AG}. This is 
true for Lie groups that arise as 
finite cartesian products of the groups $SU(r)$, $U(q)$ and $Sp(k)$ 
by \cite[Theorem 5.2]{Borel}
and thus $\Hom(\Z^{n},G)$ is path--connected for every 
$n\ge 0$. The same is true for their 
corresponding complexifications $SL_{r}(\C)$, $GL_{q}(\C)$ 
and $Sp_{k}(\C)$.  Thus $\BC G=\BC G_{\BONE}$ 
for such groups. Note that the argument provided in Theorem 
\ref{reduction to compact} works exactly in the same way if 
we replace $\BC G$ by $\BC G_{\BONE}$. Thus if 
$G$ is a real or complex reductive algebraic group and 
$K\subset G$ is a maximal compact subgroup then 
$B_{com}K_{\BONE}$  is homotopy equivalent to $\BC G_{\BONE}$.
On the other hand, define 
$E_{com}G_{\BONE}:=p^{-1}(\BC G_{\BONE})$.  Note that 
$E_{com}G_{\BONE}$ is the geometric realization of the 
simplicial subspace  of $[E_{com}G]_*$ defined by 
$[E_{com}G_{\BONE}]_n=\Hom(\Z^{n},G)_{\BONE}\times G$.
We have a 
commutative diagram 
\begin{equation*}
\begin{CD}
E_{com}G_{\BONE}@> >> EG\\
@V{p}VV     @VV{p}V\\
\BC G_{\BONE}@> >>
BG.\\
\end{CD}
\end{equation*} 
Here the lower horizontal map is the inclusion 
map $i:\BC G_{\BONE}\to BG$. After replacing $i$ 
with a fibration we obtain a fibration sequence 
$E_{com}G_{\BONE}\to \BC G_{\BONE}\to BG$
and in the same way as was done in 
Theorem \ref{reduction to compact}, we can prove that if 
$G$ is a real or complex reductive algebraic group and 
$K\subset G$ is a maximal compact subgroup then 
$E_{com}K_{\BONE}$  is homotopy equivalent to 
$E_{com}G_{\BONE}$.

If $G$ is a connected topological group, 
the long exact sequence in homotopy groups 
associated to the fibration sequence
$G\to EG\to BG$ can be used to show that
$BG$ is simply connected. Moreover, if $G$ is a simply 
connected  Lie group, then $G$ is $2$-connected since 
$\pi_{2}(G)=0$ for any Lie group $G$ (see \cite{BD}, page
225 for the compact case). Thus  $BG$ 
is $3$-connected for any such group. 
As a consequence of  \cite[Theorem 1.1]{GPS} similar 
statements are also true for $\BC G$
and $\BC G_{\BONE}$ as is proved next.

\begin{proposition}\label{3-connected}
Suppose that  $G$ is a real or complex reductive algebraic 
group that is connected as a topological space.
Then $\BC G$ and $\BC G_{\BONE}$ 
are simply connected. Moreover, if  $G$ is simply connected 
then $\BC G_{\BONE}$ is $3$-connected.  
\end{proposition}

\Proof 
By Theorem \ref{reduction to compact} we only need 
to prove the theorem for a compact connected Lie group.
Also, by Proposition \ref{proper simplicial}  in the appendix
for any Lie group $G$ the simplicial 
space  $[B_{com}G]_*$ is a proper simplicial space, in fact 
it is a strictly proper simplicial space (see Remark \ref{strictly proper}).
The same is true for $[B_{com}G_{\BONE}]_*$. 
If  $G$ is a connected Lie group 
then $[B_{com}G]_0=[B_{com}G_{\BONE}]_0=*$ is in particular 
$1$-connected and $[B_{com}G]_1=[B_{com}G_{\BONE}]_1=G$ 
is $0$-connected. By \cite[Theorem 11.12]{May} it 
follows that $\BC G$ and $\BC G_{\BONE}$ 
are simply connected. Suppose now $G$ is 
simply connected. Then $[B_{com}G_{\BONE}]_0=*$ is in particular 
$3$-connected, $[B_{com}G_{\BONE}]_1=G$ and thus this space is 
$2$-connected since $G$ is simply connected and thus 
$2$-connected as pointed out before. 
Also, $[B_{com}G_{\BONE}]_2=\Hom(\Z^{2},G)_{\BONE}$ is 
$1$-connected as it is path-connected and simply connected 
by \cite[Theorem 1.1]{GPS}. Finally, 
$[B_{com}G]_3=\Hom(\Z^{3},G)_{\BONE}$ is path-connected 
by definition. Using  \cite[Theorem 11.12]{May} it follows 
that $\BC G_{\BONE}$ is $3$-connected in this case. 
\qed

\medskip

Suppose now that $G$ is a Lie group that arises as a finite 
product of the classical groups  $SU(r)$, $U(q)$ and $Sp(k)$. 
For such groups 
$E_{com}G_{\BONE}=E_{com}G$ since $\Hom(\Z^{n},G)$ 
is path--connected for all $n\ge 0$ for such groups. Moreover, 
we have the following.
 
\begin{proposition}\label{E_{com}G 3-connected}
Assume that $G$ is a Lie group isomorphic to a finite 
product of the classical groups  $SU(r)$, $U(q)$ and $Sp(k)$ 
for $r,q,k \ge 1$.
Then $E_{com}G$ is $3$-connected.
\end{proposition}

\Proof
Observe that if $G$ and $H$ are topological groups we have 
a natural homeomorphism $\Hom(\Z^{n},G\times H) 
\cong \Hom(\Z^{n},G)\times \Hom(\Z^{n},H)$ for every 
$n\ge 0$. This implies that there is a homeomorphism 
$E_{com}(G\times H)\cong E_{com}G\times E_{com}H$. 
Because of this we only need to prove 
the proposition when $G$ is one of the groups 
$SU(r)$, $U(q)$ and $Sp(k)$. For the groups 
$SU(r)$ and $Sp(k)$ the proposition follows from the 
previous proposition and Corollary \ref{ses homotopy}. 
Thus we only need to prove the 
proposition for the case $G=U(q)$. 
By Proposition \ref{3-connected} and 
Corollary \ref{ses homotopy} it follows that $E_{com}G$ is 
simply connected. Hence to show that  
$E_{com}G$ is $3$-connected it suffices to prove that 
$\tilde{H}_{n}(E_{com}G)=0$ for $0\le n\le 3$. 
To see this, recall that 
the natural filtration of $E_{com}G$ as the geometric 
realization of the simplicial space $[E_{com}G]_{*}$ induces a 
spectral sequence 
\[
E^{2}_{p,q}=H_{p}H_{q}([E_{com}G]_{*})\Longrightarrow 
H_{p+q}(E_{com}G;\Z).
\]
The term $E^{2}_{p,q}$ in this spectral is obtained 
by taking the $p$-th homology group of the simplicial 
group $H_{q}([E_{com}G]_{*})$. Trivially we have 
$E^{2}_{0,0}=\Z$. We show next that $E^{2}_{p,q}=0$ for 
all $p,q\ge 0$ with $0<p+q\le 3$. To prove this, 
we claim that the map of simplicial spaces  
$i_{*}:[E_{com}G]_{*}\to [EG]_{*}$ induced by the inclusion 
map induces an isomorphism 
$i_{*}:H_{p}H_{q}([E_{com}G]_{*})\to H_{p}H_{q}([EG]_{*})$ 
for $0\le p+q\le 3$. Since  $H_{p}H_{q}([EG]_{*})=0$ for 
all $p+q>0$ then the proposition follows. The 
claim is trivial for $q=0$ because $[E_{com}G]_{k}$ 
is connected for all $k\ge 0$. When $q=1$ the 
simplicial groups $H_{1}([E_{com}G]_{*})$ 
and $H_{1}([EG]_{*})$ are isomorphic by 
\cite[Theorem 1.1]{GPS}.  Suppose now that $q=2$. 
Let $C_{n}=H_{2}(\Hom(\Z^{n},G)\times G;\Z)$ so that 
$\{C_{n}\}_{n\ge 0}$ is the chain complex whose $p$-th 
homology is $H_{p}H_{2}([E_{com}G]_{*})$. Trivially we have that 
$C_{0}=0$ since $H_{2}(G;\Z)=0$. Also, $C_{1}\cong \Z$ 
and the differential $\partial:C_{2}\to C_{1}\cong \Z$ is surjective 
since the inclusion $G\vee G\hookrightarrow \Hom(\Z^{2},G)$ 
induces a split injection at the level of homology by 
\cite[Theorem 1.6]{AC}. This shows that 
$H_{p}H_{2}([EG]_{*})=0$ for $p=0,1$. Finally, 
$H_{0}H_{3}([E_{com}G]_{*})$ vanishes trivially. 
\qed

\section{Commutative  K--theory}

Suppose that $X$ is a finite CW-complex and let 
$p:E\to X$ be an $n$-plane complex vector bundle.  
As in the case of a principal bundles, we say that $E$ 
is transitionally commutative if we can find an open cover 
$\{U_{i}\}_{i\in I}$ of $X$ such that 
$E$ is trivial over each $U_{i}$ and the
corresponding  transition functions commute with each other 
whenever they are simultaneously defined. This is equivalent 
to saying that, with a Hermitian metric in sight, 
the corresponding frame bundle is a 
transitionally commutative principal $U(n)$-bundle. 
Similarly, two such complex vector bundles 
$q_{0}:E_{0}\to X$ and $q_{1}:E_{1}\to X$  
are said to be transitionally commutative isomorphic if their 
frame bundles are transitionally commutative isomorphic. This 
means that we can find a transitionally commutative vector 
bundle $p:E\to X\times [0,1]$ such that 
$q_{0}=p_{|p^{-1}(X\times\{0\})}$ and $q_{1}=p_{|p^{-1}(X\times\{1\})}$.
Let $Vect_{com}(X)$ be the set of transitionally commutative 
isomorphism classes of  transitionally commutative vector bundles 
over $X$.  The Whitney sum  of two 
transitionally commutative vector bundles is also transitionally 
commutative and thus  $Vect_{com}(X)$ has the structure of 
a monoid.  We define the commutative  K--theory of $X$ to be
$K_{com}(X):=Gr(Vect_{com}(X))$, where $Gr$ denotes  
the Grothendieck construction. It is easy to see that 
if $E$ and $F$ are two  transitionally 
commutative vector bundles over $X$, then $E\oplus F$ 
is transitionally commutative isomorphic to $F\oplus E$.
This shows that,  as in the case of classical 
 K--theory, this construction defines a functor from the category of 
topological spaces to the category of abelian groups.

By Theorem \ref{geometric interpretation}  any 
transitionally commutative $n$-plane complex vector bundle 
is classified by a map $f:X\to B_{com}U(n)$. Moreover, 
two such vector bundles classified by maps $f,g:X\to B_{com}U(n)$
are transitionally commutative isomorphic 
if and only if $f$ is homotopic to $g$.  
Let $U=\colim_{n\to \infty} U(n)$, where the colimit 
is taken over the natural inclusions $i_{n}:U(n)\to U(n+1)$. 
We conclude, in an analogous way to the case of  K--theory, 
that there is a natural isomorphism of groups 
$K_{com}(X)\cong [X,\Z\times B_{com}U]$.

\begin{theorem}
The space $\Z\times \BC U$ is a loop space and for any finite
$CW$--complex $X$ there is a natural
isomorphism of groups $K_{com}(X)\cong [X, \Z\times \BC U]$.
\end{theorem}

\Proof As pointed out above we have a natural 
isomorphism $K_{com}(X)\cong [X,\Z\times B_{com}U]$  
for any finite $CW$-complex $X$. 
Consider  $M:=\bigsqcup_{n\ge 0} B_{com}U(n)$;
this space has the structure of a topological 
monoid defined as follows. For each $n,m\ge 0$ consider the 
homomorphism of topological groups
\begin{align*}
\iota_{n,m}:U(n)\times U(m)&\to U(n+m)\\
(A,B)&\mapsto \begin{bmatrix}
    A& 0\\
   0 &B 
  \end{bmatrix}.
\end{align*}
This homomorphism induces a continuous map 
\[
\Gamma_{n,m}:B_{com}U(n)\times B_{com}U(m)
=B_{com}(U(n)\times U(m))\to B_{com}U(n+m).
\]  
The different maps $\{\Gamma_{n,m}\}_{n,m\ge 0}$ 
can be assembled to obtain a  map 
$\Gamma:M\times M\to M$ giving $M$ the structure of a 
strictly associative topological monoid. Moreover, this monoid 
is commutative up to homotopy. Indeed, for each 
$n,m\ge 0$ fix a continuous path $\beta_{n,m}:[0,1]\to U(n+m)$ 
from the identity matrix $I_{n+m}$ to the matrix 
$\begin{bmatrix}
  0  &I_{m} \\
   I_{n} & 0
 \end{bmatrix}$.
Such a  path exists because $U(n+m)$ is path--connected. 
These paths induce a continuous family of homomorphisms 
\begin{align*}
h_{n,m}(t):U(n)\times U(m)&\to U(n+m)\\
(A,B)&\mapsto \beta_{n,m}(t)\iota_{n,m}(A,B)\beta_{n,m}(t)^{-1}
\end{align*}
defined for $0\le t\le 1$. After applying the functor $B_{com}$, 
these maps induce a homotopy $h:M\times M\times I\to M$ 
such that $h(A,B,0)=\Gamma(A,B)$ and $h(A,B,1)=\Gamma(B,A)$. 
The above proves that  $M$ is a strict topological monoid that is 
commutative up to homotopy. On the other hand, observe that  
$\pi_{0}(M)=\N$. Fix an element $m\in B_{com}U(1)$  
and consider the mapping telescope 
$$M_{\infty}=\text{Tel}(M\stackrel{*m}
{\rightarrow} M\stackrel{*m}
{\rightarrow} M\stackrel{*m}
{\rightarrow} \cdots)\cong  \Z\times B_{com}U.$$
By the group completion theorem (see for example \cite[Proposition 1]{McDuff}), 
it follows that the natural map 
$M\to \Omega BM$ induces a map 
$\eta:M_{\infty}\cong \Z\times B_{com}U\to \Omega B M$ 
that  is an isomorphism in homology. 
Let $(\Omega BM)_{0}$ be the path--connected 
component of $\Omega BM$ containing the trivial loop. 
Then the restriction of $\eta$, $\eta_{0}:B_{com}U\to (\Omega BM)_{0}$ 
induces an isomorphism in homology with integer coefficients.  
By Proposition \ref{3-connected}  the space 
$B_{com}U(n)$ is simply connected for every $n\ge 0$. 
Since $B_{com}U=\colim_{n\ge 0}B_{com}U(n)$ the same is true 
for $B_{com}U$.  Similarly  $(\Omega B M)_{0}$ is simply 
connected. Therefore $\eta_{0}: B_{com}U\to (\Omega BM)_{0}$ is a 
homology isomorphism between simply connected spaces. By the 
Hurewicz theorem we conclude that $\eta_{0}$ is a homotopy 
equivalence. On the other hand, since $\Omega BM$ is a 
loop space, all of its connected component are homotopy equivalent. 
We conclude that there is a homotopy equivalence 
$\Z\times B_{com}U\simeq \Omega B M$ and thus 
$\Z\times B_{com}U$ is a loop space.
\qed

\medskip

\begin{proposition}
If $X$ is a connected finite CW-complex, 
there is a natural isomorphism of groups 
$K_{com}(\Sigma X)\cong K^{0}(\Sigma X)
\times[\Sigma X,E_{com}U]$.
\end{proposition}

\Proof
By \cite[Theorem 6.3]{ACT}, for any 
$n\ge 0$ there is a natural homotopy 
equivalence given by $\theta(U(n)):U(n)\times \Omega E_{com}U(n) 
\stackrel{\simeq}{\rightarrow}  \Omega B_{com} U(n)$. 
Since $B_{com}U=\colim_{n\to \infty} B_{com}U(n)$ 
and $E_{com}U=\colim_{n\to \infty} E_{com}U(n)$ 
by passing to the colimit as $n\to \infty$ we obtain 
a homotopy equivalence 
$\Omega B_{com}U\simeq U\times \Omega E_{com}U$. 
Using this homotopy equivalence and adjunction between 
the functors $\Sigma$ and $\Omega$, we obtain 
natural isomorphisms 
\begin{align*}
K_{com}(\Sigma X)&\cong
[\Sigma X,\Z\times B_{com}U]\cong 
\Z\times [\Sigma X,B_{com}U]\\
&\cong \Z\times [X, U\times \Omega E_{com}U]
\cong \Z\times [X,U]\times [X,\Omega E_{com}U]\\
&\cong [\Sigma X, \Z\times BU]\times  [X,\Omega E_{com}U]
\cong K^{0}(\Sigma X)\times [\Sigma X, E_{com}U].
\end{align*}
\qed

\begin{example}\label{example spheres}
Using this proposition we see that  
$K_{com}(\SS^{m})\cong K^{0}(\SS^{m})$ for $0\le m\le 3$. 
For $m=0$ this is trivial; for $1\le m\le 3$ 
by the above computation there is an isomorphism 
$K_{com}(\SS^{m})\cong K^{0}(\SS^{m})
\times [\SS^{m},E_{com}U]$. 
The space $E_{com}U(n)$ is $3$-connected by 
Proposition \ref{E_{com}G 3-connected} for all $n\ge 0$. 
By passing to the colimit as $n\to \infty$ it follows that 
the same is true for $E_{com}U$. Therefore for 
$1\le m\le 3$ we have that 
$\pi_{m}(E_{com}U)\cong [\SS^{m},E_{com}U]$ 
is trivial and we conclude that 
$K_{com}(\SS^{m})\cong K^{0}(\SS^{m})$ for $0\le m\le 3$. 
However,  $K_{com}(\SS^{4})\ncong K^{0}(\SS^{4})$. 
To see this  note that the cohomological 
computations derived in Section 8 
imply that $H^{4}(E_{com}U;\Q)\ne 0$. This together with 
the universal coefficient theorem and the Hurewicz theorem 
imply that $\pi_{4}(E_{com}U)\cong [\SS^{4},E_{com}U]\ne 0$ 
and thus  $K_{com}(\SS^{4})\ncong K^{0}(\SS^{4})$. 
This in particular shows that the functor $K_{com}$ does 
not satisfy Bott periodicity for its values on spheres.
We should also mention that the non--trivial element in 
$\pi_4(B_{com}SU(2))$ which
we described in Example \ref{nontrivial class SU(2)},
is mapped non--trivially
to $\pi_4(B_{com}U(2))$ and this defines a non--trivial commutative
vector bundle over $\SS^4$ which is trivial as an ordinary bundle. 
\end{example}

\begin{remark}
In \cite{AGLT} it is proved that $\Z\times B_{com}U$ is in fact an infinite loop 
space. In particular it follows that the definition of commutative  K--theory can 
be extended to obtain a generalized cohomology theory. Moreover, it is
shown there that 
the fibration sequence 
$E_{com}U\to B_{com}U\to BU$
splits and that  
 $\Z\times B_{com}U\simeq (\Z\times BU)\times E_{com}U$ 
 as infinite loop spaces. This implies 
in particular that commutative  K--theory contains  topological  K--theory
as a summand.  Note however that, 
as proved in  Example \ref{example spheres},  commutative  K--theory 
is not $2$-periodic unlike classical   K--theory. The homotopy 
type of $\Z\times B_{com}U$ remains to be determined.
\end{remark}

\section{The topological poset associated to maximal tori in a Lie group}
\label{topological posets}

Our next goal is to provide a description of the spaces 
$B_{com}G_{\BONE}$ as a suitable homotopy colimit for any 
Lie group $G$ that is compact and connected.  To achieve this 
decomposition we first need to study the 
poset generated by all maximal tori in a compact Lie group.

\medskip

We begin by recalling the definition of a topological poset. 

\begin{definition}
A topological poset  is a partially ordered 
set $(\Ro,\preceq)$ together with a topology on the set of 
objects $\Ro$ in such a way 
that the order space
$\O_{\Ro}:=\{(x,y)\in\Ro\times \Ro ~|~ x\preceq y \}$
is a closed subset of $\Ro\times \Ro$. 
\end{definition}

A topological poset can be seen as a topological category 
where the space of objects is 
$\Ro$ and the space of morphisms is the order space $\O_{\Ro}$.  
In this article by a  topological category we mean
a small category $\Do$ for which 
the sets $\Ob(\Do)$ and $\Mor(\Do)$ 
come equipped with topologies in such a way that the 
structural maps source, target, composition and identity 
are continuous maps.

\begin{example}
Let $n\ge 0$ be a fixed integer. Given $0\le k\le n$, 
denote by $G_{k}(\C^{n})$ the Grassmannian manifold consisting 
of all those $k$-dimensional $\C$-vector spaces in $\C^{n}$. 
Let $\Gr(n)$ be the poset of all $\C$-vector subspaces 
in $\C^{n}$. This set naturally has the structure of a poset 
by inclusion. Note that 
$\Gr(n)=\bigsqcup_{0\le k\le n}G_{k}(\C^{n})$.
We can use this identification to give $\Gr(n)$ a topology making 
it into a topological poset.
\end{example}

The maximal tori in a compact Lie group
$G$ define a topological poset in the following natural way.

\begin{definition}\label{poset T(G)}
Suppose that $G$ is a Lie group. Define $\T(G)$ to be the
poset whose objects are closed subspaces 
$S\subset G$ arising 
as the intersection of a collection of maximal tori in $G$, 
with the order 
in $\T(G)$ given by inclusion.
\end{definition} 

The set  $\T(G)$ can be given a topology making 
it into a topological poset as follows.
Let $\Co(G)$ denote the set  of all  closed subspaces in 
$G$.   Suppose that  $\U:=\{U_{1},\dots,U_{n}\}$ 
is a finite collection of  open sets in $G$.  Define 
$\Co(G,\U)$ to be the set of elements $A\in \Co(G)$ such that 
$A\subset \cup_{i=1}^{n}U_{i}$ and $A\cap U_{i}\ne \emptyset$ for 
all $1\le i\le n$. The sets of the form $\Co(G,\U)$ form a basis for a 
topology in $\Co(G)$ called the finite topology 
(see \cite{Michael} for details). Note that $\T(G)\subset \Co(G)$ 
and in this way we can give $\T(G)$ the subspace topology making 
it into a topological poset.  Our next goal is to describe
the structure of $\T(G)$ as a topological 
space for any compact connected Lie group $G$. 

Let $\g$ denote the Lie algebra 
of $G$ and  fix $T\subset G$ a maximal torus in $G$ with Lie algebra 
$\t$. Let $\Phi$ be the root system 
associated to $(\g,\t)$ and fix a subset $\Phi^{+}$ of positive 
roots of $\Phi$.  For each $\alpha\in \Phi^{+}$ and any integer $n$, 
define 
$$\t_{\alpha}:= \{ X\in \t ~|~ \alpha(X)\in 2\pi i \Z \},\,\,\,
\t_{\alpha,n}:= \{ X\in \t ~|~ \alpha(X)=2\pi i n \}.$$
Each $\t_{\alpha,n}$ is a hyperplane of codimension one and the set 
$D(G):=\displaystyle{\bigcup_{\alpha \in \Phi^{+}}\t_{\alpha}}$
is  called the Stiefel diagram of $G$. Recall that an element $g\in G$ 
is called singular if it belongs to more than one maximal torus in $G$. 
Equivalently, $g\in G$ is singular if and only if 
${\rm{dim}}(Z_{G}(g))>{\rm{dim}}(T)$.  
Let $G_{s}\subset G$ be the set of  singular elements in $G$ 
and let $T_{s}=T\cap G_{s}$ be the set of singular elements in 
$T$. Consider the restriction of the exponential map $\exp:\t \to T$.  
We have $\exp^{-1}(T_{s})=D(G)$; that is,  if $X\in \t$ then 
$\exp(X)$ is singular if and only if $X\in D(G)$.  
Given a set of positive roots 
$I=\{\alpha_{1},\dots, \alpha_{k}\}$ define 
\[
t_{I}=\bigcap_{i=1}^{k}t_{\alpha_{i}} \text{ and }
T_{I}:=\exp(t_{I})\subset T. 
\]
In the previous definition we allow the case $k=0$. In this case 
we take the convention that  $t_{\emptyset}=\t$ and thus 
$T_{\emptyset}=T$. Let $\Delta=\{\alpha_{1},\dots ,\alpha_{r}\}$ be a 
set of simple roots for the root system $\Phi$.  Recall that the 
Weyl group $W$ is a reflection subgroup generated by the 
reflections $s_{\alpha}$ corresponding to  elements 
$\alpha \in \Delta$.  Given $I\subset \Delta$, the subgroup of 
$W_{I}$ of $W$ generated by the reflections $s_{\alpha}$ 
corresponding to elements $\alpha\in I$ is called a parabolic 
subgroup of $W$. Note that each parabolic subgroup is 
itself a reflection subgroup. We are interested in the different 
parabolic subgroups of $W$ up to conjugacy.  If 
$I,J\subset \Delta,$ then $W_{I}$ is conjugated to $W_{J}$ 
if and only if $I=wJ$ for some $w\in W$.  In that case we 
say that  $I$ and $J$ are in the 
same Coxeter class and  write $I\sim_{W} J$. The 
relation $\sim_{W}$ defines an equivalence relation on the 
set of subsets of $\Delta$. We denote by $\E_{W}$ the set of 
equivalence classes of subsets of $\Delta$ under this equivalence 
relation and by $[I]$  the equivalence class in $\E_{W}$ that 
contains $I\subset \Delta$.

\begin{theorem}\label{structure T(G)}
Suppose that $G$ is a compact connected Lie group.
Fix $\Delta=\{\alpha_{1},\dots ,\alpha_{r}\}$ a set of simple roots. 
Then any element $S\in \T(G)$ is conjugated to  
$T_{I}$ for some  
$I\subset \Delta$. Moreover, there is a  
$G$-equivariant homeomorphism 
$\T(G)\cong\bigsqcup_{[I]\in \E_{W}}
G/N_{G}(T_{I})$. 
\end{theorem}

\Proof Fix a maximal torus $T\subset G$ and suppose that 
$S\in \T(G)$. Since any two maximal tori in $G$ are conjugated,  
after replacing $S$ with  a suitable conjugate we may 
assume that $S\subset T$. If $S=T$, then $S=T_{\emptyset}$ 
and there is nothing to prove. 
Assume then that $S\subsetneq  T$. Let's show that under this 
assumption $S$ is conjugated to $T_{I}$ for some set of simple roots 
$I\subset \Delta$.   Let $S_{0}$ be the identity component of $S$. 
Then $S_{0}$ is a compact, connected and abelian subgroup of $G$. 
Thus $S_{0}$ is a torus and in particular we can find an element 
$x_{0}\in S_{0}$ such that $S_{0}=\overline{\left< x_{0}\right>}$. Let 
$\beta_{1},\dots, \beta_{l}$ be the set of positive roots $\alpha$ with 
$x_{0}\in T_{\alpha}$. It follows that 
$S_{0}\subset  T_{\beta_{1}}\cap\cdots\cap T_{\beta_{l}}=T_{J}$, 
where $J=\{\beta_{1},\dots,\beta_{l} \}$. In fact $S_{0}\subset T_{J,0}$, 
where $T_{J,0}$ denotes the identity component of $T_{J}$.  
As a first step we show that $S_{0}=T_{J,0}$. To see this recall that 
the adjoint representation provides a decomposition 
of the complexification of $\g$ into  a direct sum of weight spaces
$\g_{\C}=\t_{\C}\oplus \left( \bigoplus_{\alpha\in \Phi} L_{\alpha}\right)$.
This in turn provides a decomposition of $\g$  in the form 
$\g=\t\oplus \left( \bigoplus_{\alpha\in \Phi^{+}} M_{\alpha}\right),$
where $M_{\alpha}=(L_{\alpha}\oplus L_{-\alpha})\cap \g$. 
Note that $Z_{G}(x_{0})=Z_{G}(S_{0})$ and in particular this group 
is connected since the centralizer of any torus in $G$ is connected. 
By  \cite[Proposition V 2.3]{BD} the Lie algebra of $Z_{G}(x_{0})$ is 
$\z_{\g}(x_{0})=\t \oplus \left(\bigoplus_{i=1}^{l} M_{\beta_{i}}\right)$. 
On the other hand, the Lie algebra of $T_{J,0}$ is 
$\t_{J,0}:= \t_{\beta_{1},0}\cap\cdots\cap \t_{\beta_{l},0}$. It follows that 
the Lie algebra of $Z_{G}(T_{J,0})$ is 
$\z_{\g}(\t_{J,0})=\t \oplus \left(\bigoplus_{i=1}^{l} M_{\beta_{i}}\right)$. 
This proves that  $Z_{G}(x_{0})$ and $Z_{G}(T_{J,0})$ are connected 
subgroups of $G$ that have the same Lie algebra which in turn implies 
that $Z_{G}(x_{0})=Z_{G}(T_{J,0})$.  Suppose 
now that $T'$ is a maximal torus that contains 
$x_{0}$. Then $T'\subset Z_{G}(x_{0})=Z_{G}(T_{J,0})$ which implies 
that $T_{J,0}\subset T'$ since the centralizer of a connected abelian 
subgroup of $G$ is the union of all maximal tori containing it.
This shows that $T_{J,0}$ is contained in the intersection of al maximal 
tori that contain $x_{0}$. Since $S$ is the intersection of a family of 
maximal tori, it follows that $T_{J,0}\subset S$ and by 
connectedness we have $T_{J,0}\subset S_{0}$. We conclude 
that $S_{0}=T_{J,0}$. To show that $T_{J}= S$ 
recall that the center of $G$ is the intersection of all maximal tori in $G$. 
Therefore $S$  contains the center of $G$. This is also true for $T_{J}$. 
Moreover, it is easy to see that the center of $G$ intersects all the 
connected components of $S$ and $T_{J}$. From here it follows that 
$S=T_{J}$ and that the Lie algebra of $S$ is $\t_{J,0}$.
Let's show now that after replacing $S$ with some further conjugate, 
we can  choose the $\beta_{j}$'s to be simple roots. 
Choose a minimal set of positive roots $\gamma_{1},\dots, \gamma_{k}$ 
with $\t_{J,0}=\t_{\gamma_{1},0}\cap\cdots \cap \t_{\gamma_{k},0}$. 
Then we have proper inclusions 
$\t_{\gamma_{1},0}\cap\cdots \cap \t_{\gamma_{k},0}\subset 
\t_{\gamma_{1},0}\cap\cdots \cap \t_{\gamma_{k-1},0}
\subset\cdots\subset \t_{\gamma_{1,0}}$.
As $\t\setminus \bigcup_{\alpha\in \Phi^{+}}\t_{\alpha,0}$ is a union 
of Weyl chambers, then we can find some (closed) Weyl chamber 
$\mathfrak{C}$ in such a way that each 
$\t_{\gamma_{i},0}$ is a face of 
$\mathfrak{C}$ for every $1\le i\le k$. Associated to the Weyl chamber  
$\mathfrak{C}$ there is a set of simple roots of $\Phi$. Since each 
$\t_{\gamma_{i},0}$ is a face of the chamber $\mathfrak{C}$ then, 
after replacing the sign of the $\gamma_{i}$'s if necessary, 
the roots  $\gamma_{1},\dots, \gamma_{k}$  are roots in some base of 
$\Phi$. Since the Weyl group acts transitively on the set of all bases 
on $\Phi$, it follows that  we can find some $w\in W$ such that  
$w\gamma_{1},\dots, w\gamma_{k}$  are  in $\Delta$. This shows that 
$S$ is conjugated to 
$T_{I}$, where $I
=\{\alpha_{i_{1}},\dots, \alpha_{i_{k}}\}$ is some set of 
of simple roots.

On the other hand, it is easy to see that for any 
$I\subset \Delta$, then 
the closed subspace 
$T_{I}$ is the intersection 
of all the maximal tori containing it and thus  
$T_{I}\in \T(G)$.  Also, 
if $I$ and $J$ are subsets of $\Delta$, 
then $T_{I}$  is conjugated to $T_{J}$ 
if and only if $J=wI$ for some 
$w\in W$; that is, $T_{I}$  is conjugated to 
$T_{J}$  if and only if $I$ and 
$J$ are in the same Coxeter class.
To finish 
note that the space of subgroups in $G$ that are conjugated to 
$T_{I}$  can be identified with 
$G/N_{G}(T_{I})$ and  
thus the theorem follows.
\qed

\medskip

\begin{example}\label{poset case U(n)} Suppose that $G=U(n)$ 
for $n\ge 1$. For this group a maximal torus $T$ can be chosen 
to be the set of all diagonal matrices with diagonal entries 
$x_{1},\dots, x_{n}\in \SS^{1}$. The Weyl group $W=\Sigma_{n}$ 
acts by permutation on the diagonal entries.
The Lie algebra $\t$ can be identified with 
$\t=\{(X_{1},\dots,X_{n}) ~|~  X_{i}\in i\R 
\text{ for all } 1\le i\le n \}$.
The root system $\Phi$ consists of all functions 
$\alpha_{i,j}(X_{1},\dots,X_{n}) =X_{i}-X_{j}$
for $1\le i, j\le n$ with $i\ne j$. 
A choice of positive roots is the set of roots $\alpha_{i,j}$ with 
$i<j$ and the roots 
\[
\Delta:=\{\alpha_{1}:=\alpha_{1,2},  \alpha_{2}:=\alpha_{2,3}, \dots,  
\alpha_{n-1}:=\alpha_{n-1,n}\}
\]
form a set of simple roots. By the previous 
theorem any $S\in \T(U(n))$ is conjugated to some 
$T_{I}$, where $I:=\{\alpha_{i_{1}},
\dots,  \alpha_{i_{k}}\}$ and
$1\le i_{1}<\cdots <i_{k}\le n-1$. Unraveling the definition, we see 
that $T_{I}$ consists of all 
diagonal matrices with entries $x_{1},\dots, x_{n}$ 
and with $x_{i_{r}}=x_{i_{r}+1}$ for all $1\le r\le k$. In other words, 
the roots  $\alpha_{i_{1}},\dots,  \alpha_{i_{k}}$ determine the number 
of repeated diagonal entries in the elements of 
$T_{I}$.  The conjugacy classes of such tori can be parametrized 
using partitions of the number $n$.
Recall that a nondecreasing sequence of integers 
$\lambda=(\lambda_{1},\dots, \lambda_{k})$ is a partition of $n$ if 
$n=\lambda_{1}+\lambda_{2}+\cdots+\lambda_{k}$. We write
$\lambda \vdash n$ to mean that $\lambda$ is a partition of $n$. 
The set of Coxeter classes in 
$\Delta$ are in bijective correspondence with the set of partitions  
$\lambda$ of the number $n$. This can be seen by corresponding 
for a partition $\lambda$ of $n$ the Coxeter class represented by  
the set of simple roots
$I(\lambda):=\Delta\setminus\{\alpha_{\lambda_{1}},\dots, 
\alpha_{\lambda_{1}+\cdots+\lambda_{k-1}} \}$.
Note that if $\lambda \vdash n$ then the 
associated torus $T_{I(\lambda)}$ is the subspace 
of $T$ consisting of those diagonal matrices with entries of the 
form 
\[
(\underbrace{x_{1},\dots,x_{1}}_{\lambda_{1}},\dots,
\underbrace{x_{k},\dots,x_{k}}_{\lambda_{k}}).
\] 
Therefore each $S\in \T(U(n))$ is conjugated to 
$T_{I(\lambda)}$ for some unique $\lambda \vdash n$.
On the other hand,  given a partition $\lambda$ of $n$, let 
$Fl(\lambda):=U(n)/(U(\lambda_{1})\times\cdots 
\times U(\lambda_{k}))$. The space  
$Fl(\lambda)$ is the  flag manifold consisting of all flags 
of the form 
$V_{1}\subset V_{2}\subset\cdots \subset V_{k}=\C^{n}$,
where $V_{i}$ is a $\C$-vector subspace of $\C^{n}$ of 
dimension 
${\rm{dim}}_{\C}(V_{i})=\lambda_{1}+\cdots+\lambda_{i}$. 
We can see the partition $\lambda$ as an ordered 
$k$-tuple $(\lambda_{1},\dots,\lambda_{k})$. The symmetric 
group $\Sigma_{k}$ acts on the set of such $k$-tuples 
by permutation. We denote by $\Sigma_{\lambda}$ 
the isotropy of $\Sigma_{k}$ at $\lambda$ under this action. 
With this notation, the group $N_{U(n)}(T_{I(\lambda)})$ 
fits in a short exact sequence
\[
1\to U(\lambda_{1})\times\cdots 
\times U(\lambda_{k}) \to N_{U(n)}(T_{I(\lambda)})\to 
\Sigma_{\lambda}\to 1.
\]
It follows that if $\lambda \vdash n$, then 
\[
U(n)/N_{U(n)}(T_{I(\lambda)})\cong(U(n)/(U(\lambda_{1})\times\cdots 
\times U(\lambda_{r})))/\Sigma_{\lambda}=Fl(\lambda)/\Sigma_{\lambda}.
\]
We conclude that 
$\T(U(n))=\bigsqcup_{\lambda\vdash n}
Fl(\lambda)/\Sigma_{\lambda}$, 
where $\lambda$ runs through all partitions of $n$.
\end{example}

\section{The space $B_{com}G_{\BONE}$ as a homotopy colimit}
\label{homotopy colimits} 

In this section we derive a description of $B_{com}G_{\BONE}$ 
as a suitable homotopy colomit for a real or complex reductive algebraic group $G$.  
Note that by Theorem \ref{reduction to compact} we can work 
without loss of generality in the category of compact Lie groups. 

\medskip

To start we show that the space $B_{com}G_{\BONE}$ is a 
colimit over the poset $\T(G)$.

\begin{proposition}\label{property max tori}
For any compact connected Lie group $G$ we have 
$B_{com}G_{\BONE} \cong\colim_{S\in \T(G)}BS$.
\end{proposition}

\Proof Suppose that $T\subset G$ is a maximal torus. Then 
$T^{n}\subset \Hom(\Z^{n},G)_{\BONE}$ for all $n\ge 0$ 
and thus $BT\subset B_{com}G_{\BONE}$. This proves that 
$\bigcup_{T\in \T(G)} BT\subset B_{com}G_{\BONE}$.
Suppose now that $x\in  B_{com}G_{\BONE}$. 
Then $x$ can be represented 
in the form $x=[(g_{1},\dots,g_{n},t)]$ for some 
$(g_{1},\dots,g_{n})\in \Hom(\Z^{n},G)_{\BONE}$ and $t\in  \Delta_{n}$. 
By \cite[Lemma 4.2]{Baird}  there is a maximal torus 
$T\subset G$ such that $g_{i}\in T$ for all $1\le i\le n$. 
Therefore  $x\in BT$.  This proves that 
$B_{com}G_{\BONE} =\bigcup_{T\in\T(G)} BT$ and  
thus  
$B_{com}G_{\BONE}= 
\bigcup_{T\in\T(G)} BT\cong \text{colim}_{S\in \T(G)}BS$.
\qed

\medskip

As is well known homotopy colimits are better suited for 
homotopical computations than colimits. Therefore we would 
like to obtain a decomposition $B_{com}G_{\BONE}$ as a homotopy 
colimit over a suitable category. It can be seen that 
the space $B_{com}G_{\BONE}$ can be described 
as the homotopy colimit of the functor $B:\T(G)\to \Top$. 
However, this decomposition is not very helpful as the category 
$\T(G)$ is a topological category (see \cite{Hollender} and \cite{Lind} for 
a discussion of homotopy colimits over topological categories). 
In particular the usual Bousfield-Kan 
spectral sequence does not apply to such homotopy colimits. 
We will get around this issue by obtaining a decomposition of 
$B_{com}G_{\BONE}$ as a homotopy colimit over a 
\textit{discrete} category.  A key element in this decomposition 
is the rank function defined over the poset $\T(G)$. To be more 
precise, suppose that $G$ is a compact connected Lie group. Let 
$Z=Z(G)$ be the center of $G$ and write 
$n={\rm{rank}}(G)-{\rm{rank}}(Z)\ge 0$. If $S\in \T(G)$ 
then $S$ is a closed subgroup of $G$ 
and in particular it is a compact Lie group. Define 
\[
\rho(S):={\rm{rank}}(S)-{\rm{rank}}(Z).
\]
This way for every $S\in \T(G)$ we obtain  $0\le \rho(S)\le n$.

\begin{proposition}\label{function rho}
If $G$ is a compact connected Lie group the  
function $\rho:\T(G)\to \N$  is strictly increasing and  
continuous. Moreover, $\rho$ attains any value  $0\le m\le n$.
\end{proposition}

\Proof 
Suppose that  $S_{2}\subsetneq S_{1}$ are elements 
in $\T(G)$. Fix a maximal torus $T\subset G$ 
and  a set $\Phi^{+}$ of positive roots. 
By Theorem \ref{structure T(G)} 
we can find simple roots 
$I:=\{\alpha_{i_{1}},\dots, \alpha_{i_{k}}\}$ 
and ${J}:=\{\beta_{j_{1}},\dots, \beta_{j_{l}}\}$
and $g_{1},g_{2}\in G$ such that 
$S_{2}=g_{2}T_{J}g_{2}^{-1}\subsetneq 
S_{1}=g_{1}T_{I}g_{1}^{-1}$.
Since the exponential map is surjective this implies  
$g_{2}\t_{J}g_{2}^{-1}
\subsetneq g_{1}\t_{I}g_{1}^{-1}$.
In fact we must have $g_{2}\t_{J,0}g_{2}^{-1}
\subsetneq g_{1}\t_{I,0}g_{1}^{-1}$. 
Therefore ${\rm{rank}}(S_{2})
={\rm{dim}}_{\R}(g_{2}\t_{J,0}g_{2}^{-1})
<{\rm{dim}}_{\R}(g_{1}\t_{I,0}g_{1}^{-1})={\rm{rank}}(S_{1})$ 
and thus $\rho(S_{2})<\rho(S_{1})$. Also recall 
that by Theorem \ref{structure T(G)} we have 
$\T(G)\cong\bigsqcup_{[I]\in \E_{W}}G/N_{G}(T_{I})$. 
The map $\rho$  is constant on each 
connected component of the form $G/N_{G}(T_{I})$. 
Thus $\rho:\T(G)\to \N$ is a continuous function. 
Finally, fix $0\le m\le n$ and let $r$ be the rank 
of the center of $G$. Choose any set $I$ consisting of  
$m+r$ simple roots. Then the corresponding element 
$T_{I}\in \T(G)$ is such that $\rho(T_{I})=m+r-r=m$.
\qed

\medskip

Choose $n\ge 0$  as in the previous proposition and let  
$\So(n)$  be the poset 
consisting of all the nonempty subsets of $\{0,1,\dots, n\}$, with 
the order given by the \textit{reverse} inclusion of sets. 
We see an element in 
$\So(n)$ in the form $\i:=\{i_{0},\dots,i_{k}\}$, 
with $0\le i_{0}<i_{1}<\cdots<i_{k}\le n$.  
Associated to the group $G$ we 
have a functor 
$\Fo_{G}:\So(n)\to \Top$ defined in the following way. 
Suppose that 
$\i:=\{i_{0},\dots,i_{k}\}$ is an object in $\So(n)$. Then 
we define 
\[
\Fo_{G}(\i):=\{(S_{0},\dots, S_{k},a) 
~|~ S_{0}\subset\cdots\subset S_{k}\in \T(G), \  
\rho(S_{r})=i_{r} \text{ for } 0\le r\le k 
\text{ and } a\in BS_{0} \}. 
\]
Note that we have a natural inclusion 
$\Fo_{G}(\i)\subset \T(G)^{k+1}\times BG$ 
and we give $\Fo_{G}(\i)$ the subspace topology. 
If $\j$ is a subset of $\i$ then the 
natural projection maps induce continuous functions 
$p_{\i,\j}:\Fo_{G}(\i)\to \Fo_{G}(\j)$. This defines a 
functor $\Fo_{G}:\So(n)\to \Top$.   To simplify matters, 
we use the following notation for elements in $\Fo_{G}(\i)$. 
Given a chain $S_{0}\subset\cdots\subset S_{k}$  in $\T(G)$ 
with $\rho(S_{k})=i_{r}$ for  $0\le r\le k$, we denote 
$\Si_{\i}=(S_{0},\dots, S_{k})$. With this notation the 
objects in $\Fo_{G}(\i)$ are pairs of the form 
$(\Si_{\i},a)$ with $a\in BS_{0}$.

\begin{theorem}\label{replacing discrte homotopy colims}
Suppose that $G$ is a compact connected Lie group. 
Then there is a natural homotopy equivalence
\[
\hocolim_{\i\in \So(n)}\Fo_{G}(\i)\simeq  B_{com}G_{\BONE}.
\]
\end{theorem}
\Proof
The proof of this theorem will be divided into two steps. 
As a first step we construct a topological category 
$\Do$ in such a way that there is a homotopy equivalence 
$B\Do \simeq B_{com}G_{\BONE}$. 

Let $\Co$ be the topological category  whose objects are 
the elements in $B_{com}G_{\BONE}$ and the only 
morphisms are the identity maps. 
Since there are 
no nontrivial morphisms in $\Co$ we have 
$B\Co=B_{com}G_{\BONE}$ as topological spaces. 
Also, consider the topological category $\Do$ defined as 
follows. The objects in $\Do$ are pairs of the form 
$(S,a)$, where $S\in \T(G)$  and $a\in BS$. If 
$(S_{1},a)$ and $(S_{2},b)$ are two objects in $\Do$, 
then there is a unique morphism $i:(S_{1},a)\to (S_{2},b)$ 
if and only if $a=b$ and $S_{1}\subset S_{2}$.  
We have a functor $F:\Do\to \Co$ that for  an object 
$(S,a)$ in $\Do$ corresponds 
$F(S,a)=a\in BS\subset B_{com}G_{\BONE}$.  
The functor $F$ sends every morphism in $\Do$ to the 
correspondent identity morphism in $\Co$. Fix an element  
$a\in B_{com}G_{\BONE}$ and consider the under category 
$a\backslash F$. The objects in this 
category are tuples of the form $((S,a),\text{id}_{a})$, 
where $S\in \T(G)$ is such that $a\in BS$ and $\text{id}_{a}:a\to a$ is the 
identity morphism in $\Co$. There is a unique 
morphism $((S_{1},a),\text{id}_{a})\to ((S_{2},a),\text{id}_{a})$ 
in $a\backslash F$ whenever $S_{1}\subset S_{2}$. 
We observe that the  category $a\backslash F$ has an 
initial object. Indeed,  let $S_{a}=\bigcap_{S\in \T(G), a\in BS}S$. 
Then $S_{a}$ is the  smallest element 
in $\T(G)$ such that  $a\in BS_{a}$ and this implies that 
$((S_{a},a),\text{id}_{a})$  is an initial object in $a\backslash F$. 
We conclude that the category $a\backslash F$ is a contractible 
category since it has an initial object. This means that is classifying 
space is contractible. 
Therefore as a particular case of \cite[Lemma A.5]{Lind} 
we obtain that the map 
$BF:D\Do\to B\Co\cong B_{com}G_{\BONE}$ 
is a homotopy equivalence. 

As a second step we show that there is a natural homeomorphism 
$\phi:\hocolim_{\i\in \So(n)}\Fo_{G}(\i)\to B\Do$. This will finish the 
proof.
For this recall that by definition
$\hocolim_{\i\in \So(n)}\Fo_{G}(\i)=B\Go$, where 
$\Go$ is the topological category 
whose objects are sequences of the form 
$(\i,\Si_{\i},a)$. Here $\i=\{i_{0}\dots,i_{k}\}$ is 
an element in $\So(n)$, $\Si_{\i}=(S_{0},\dots, S_{k})$ 
is a chain in $\T(G)$ with $\rho(S_{r})=i_{r}$ for 
$0\le r\le k$ and $a\in BS_{0}$.  Whenever $\j\subset \i$ 
there is a unique morphism in $\Go$,  $(\i,\Si_{\i},a)\to (\j,\Si_{\j},a)$, 
which is induced by the corresponding projections. 
An element in $B\Go$ is of the form 
$w=[(g_{1},\dots,g_{l}),t]$, where $t\in \Delta_{l}$ and 
$(g_{1},\dots,g_{l})$ is a sequence of composable 
morphisms in $\Go$ of the form
\[
(\i_{l},\Si_{\i_{l}},a)\stackrel{g_{l}}{\longrightarrow}
(\i_{l-1},\Si_{\i_{l-1}},a)\stackrel{g_{l-1}}{\longrightarrow}
\cdots\stackrel{g_{1}}{\longrightarrow} (\i_{0},\Si_{\i_{0}},a) .
\]
 This implies in particular   that
$\i_{0}\subset \cdots\subset \i_{l}$ is a nested sequence of 
nonempty sets. Write $\i:=\i_{l}=\{i_{0},\dots,i_{k}\}$ 
and $\Si_{\i}=(S_{0},\dots, S_{k})$. 
 Therefore the 
morphisms $g_{1},\dots,g_{l}$ induce composable 
morphisms in $\Do$ 
\[
(S_{0},a)\stackrel{f_{1}}{\longrightarrow}
(S_{1},a)\stackrel{f_{2}}{\longrightarrow}
\cdots\stackrel{f_{k}}{\longrightarrow} (S_{k},a),
\]
where $f_{r}$  is the unique morphism in 
$\Do$ from $(S_{r-1},a)$ to $(S_{r},a)$. Consider 
the standard $k$-simplex $\Delta_{k}$ that 
corresponds to the composable sequence 
$(f_{k},\dots, f_{1})$ in $\Do$.  
We identify the vertices of  $\Delta_{k}$ 
with the numbers $i_{0},\dots, i_{k}$ by corresponding to each
$i_{r}$ the vertex $v_{i_{r}}=(0,\dots,1,\dots,0)\in \Delta_{k}$,  
with the entry $1$ in the $r$-th position.  
If $\j=\{i_{r_{0}},\dots, i_{r_{m}}\}\subset \i$,  we denote by  
$v_{\j}$ the barycenter of the simplex generated by the 
vertices $v_{i_{r_{0}}},\dots, v_{i_{r_{m}}}$. Thus with this notation, 
all the vertices in the barycentric subdivision 
of $\Delta_{k}$, $\Bo\Delta_{k}$,  are of the form $v_{\j}$, 
where $\j$ is a nonempty subset of $\i$. We can  
associate to the nested sequence 
$\i_{0}\subset \cdots\subset \i_{l}$ the face of $\Bo\Delta_{k}$ 
generated by the vertices $v_{\i_{0}},\dots, v_{\i_{l}}$.  Given 
$t=(t_{0},\dots,t_{l})\in \Delta_{l}$, define 
$\gamma(t)\in \Delta_{k}$ by 
\[
\gamma(t)=\gamma(t_{0},\dots, t_{l})=t_{0}v_{\i_{l}}+i_{1}v_{\i_{l-1}}+\dots+
t_{l}v_{\i_{0}}.
\]
Using this convention we define 
\[
\phi([(g_{1},\dots,g_{l}),t]):=[(f_{k},\dots, f_{1}),\gamma(t)]\in B\Do.
\]
In other words, the map $\phi$ is a linear isomorphism from the  
standard simplex $\Delta_{l}$, that corresponds to the composable 
sequence $(g_{1},\dots,g_{l})$ in $\Go$, onto the face of 
$\Bo\Delta_{k}$ that corresponds to the chain 
$\i_{0}\subset \i_{2}\subset \cdots \subset \i_{l}$,
where $\Delta_{k}$ is the simplex associated to the composable 
sequence $(i_{k},\dots, i_{1})$ in $\Do$. 
It can be seen that the map $\phi$ is well defined and 
is in fact a homeomorphism.
\qed

\medskip

The values of the functor $\Fo_{G}$ can be described explicitly 
as follows.  Fix $\i=\{i_{0},\dots,i_{k}\}$ an  element 
in $\So(n)$. The conjugation action of $G$ 
defines an equivalence relation 
on the set of chains in $\T(G)$ in the following way. 
Suppose that $\Si_{\i}=\{S_{0}\subset\cdots\subset S_{k}\}$  
and $\Si'_{\i}=\{S'_{0}\subset\cdots\subset S'_{k}\}$ are 
two chains with $\rho(S_{r})=\rho(S'_{r})=i_{r}$ for $0\le r\le k$. 
Then we say that $\Si_{\i}\sim \Si'_{\i}$ if and only if we can find some 
$g\in G$ such that $\Si_{\i}=g\Si'_{\i}g^{-1}$; that is, 
$S_{r}=gS'_{r}g^{-1}$ for all $0\le r\le k$. Denote 
by $\E(\i)$  the set of all equivalence classes of such chains 
and by $[\Si_{\i}]$ the equivalence class representing $\Si_{\i}$ 
in $\E(\i)$. 
Fix a chain $\Si_{\i}:=\{S_{0}\subset\cdots\subset S_{k}\}$ 
in $\T(G)$ with $\rho(S_{r})=i_{r}$ for $0\le r\le k$. The 
conjugation action of $G$ induces a continuous map 
\begin{align*}
\bar{\mu}_{\Si_{\i}}: G\times BS_{0}&\to \Fo_{G}(\i)\\
(g,a)&\mapsto (g\Si_{\i}g^{-1},gag^{-1}).
\end{align*}
Let  $N_{G}(\Si_{\i})$ be the normalizer of $\Si_{\i}$ in $G$; that 
is, the subgroup of $G$ consisting of elements $g\in G$ such that
$g\Si_{\i}g^{-1}=\Si_{\i}$. The group $N_{G}(\Si_{\i})$ 
acts by conjugation on $BS_{0}$ and on the left on $G$ by 
the assignment $n\cdot g =gn^{-1}$. This induces a diagonal 
action of $N_{G}(\Si_{i})$ on $G\times BS_{0}$ and the map 
$\bar{\mu}_{\Si_{\i}}$ is invariant under this action. Therefore 
$\bar{\mu}_{\Si_{\i}}$ induces  a continuous function
\[
\mu_{\Si_{\i}}:G\times_{N_{G}(\Si_{\i})} BS_{0}\to \Bo_{\T(G)}(\i).
\]
If we vary $\Si_{\i}$ through the different equivalence classes 
in $\E(\i)$, then we obtain  a continuous map 
\[
\mu_{\i}=\bigsqcup_{[\Si_{\i}]\in \E(\i)}\mu_{\Si_{\i}}:
\bigsqcup_{[\Si_{\i}]\in \E(\i)}G\times_{N_{G}(\Si_{\i})} BS_{0}
\to  \Fo_{G}(\i).
\]
This map is clearly surjective. In fact this map is also injective. 
Indeed, suppose that 
$\mu_{\i}(g,a)=(g\Si_{\i}g^{-1},gag^{-1})
=(g_{1}\Si'_{\i}g_{1}^{-1},g_{1}a_{1}g_{1}^{-1})=\mu_{\i}(g_{1},a_{1})$. 
Then  $g\Si_{\i}g^{-1}=g_{1}\Si'_{\i}g_{1}^{-1}$ which means that 
$[\Si_{\i}]= [\Si'_{\i}]$. Thus we can assume without loss of generality 
that $\Si_{\i}=\Si'_{\i}$. Also we have  
$g_{1}^{-1}g\Si_{\i}(g_{1}^{-1}g)^{-1}=\Si_{\i}$ and 
$gag^{-1}=g_{1}a_{1}g_{1}^{-1}$. Therefore 
$n:=g_{1}^{-1}g\in N_{G}(\Si_{\i})$ is such that $nan^{-1}=a_{1}$.
We conclude that in 
$G\times_{N_{G}(\Si_{\i})}BS_{0}$ we have 
$[(g,a)]=[gn^{-1},nan^{-1}]=[(g_{1},a_{1})]$, proving that  $\mu_{\i}$ 
is injective. Moreover, it can easily be seen that $\mu_{\i}^{-1}$ 
is also continuous and thus $\mu_{\i}$ is a homeomorphism.
We conclude that for every element $\i$ in $\So(n)$ 
there is a natural homeomorphism
\[
\Fo_{G}(\i)\cong \bigsqcup_{[\Si_{\i}]\in \E(\i)}G\times_{N_{G}(\Si_{\i})} BS_{0}.
\]

The sets $\E(\i)$ that appear in the previous description 
can be expressed in terms of the root system 
$\Phi$ associated to a maximal torus $T\subset G$ 
in the following way. Let 
$\i=\{i_{0},\dots,i_{k}\}$ be an object in $\So(n)$ and 
$\Si_{\i}=\{S_{0}\subset\cdots\subset S_{k}\}$  
a chain in $\T(G)$ with  $\rho(S_{r})=i_{r}$ for $0\le r\le k$.  
After replacing $\Si_{\i}$ with a suitable conjugate we may assume that 
$\Si_{\i}$ is such that 
$S_{0}\subset\cdots\subset S_{k}\subset T$. By Theorem \ref{structure T(G)} 
for every $0\le r\le k$ we can find a set of simple roots 
$I_{r}$ in such a way that $S_{r}=g_{r}T_{I_{r}}g_{r}^{-1}$ for 
some $g_{r}\in G$.  Let $T_{I_{r},0}$ denote 
the connected component of $T_{I_{r}}$ that contains the identity. 
Then $T_{I_{r},0}$ is a torus since it is a compact, connected abelian Lie group. 
Therefore for each $0\le r\le k$ we can find some element 
$x_{r}$ such that $T_{I_{r},0}=\overline{\left< x_{r}\right>}$.  
Each  $x_{r}$ is such that $x_{r}\in T_{I_{r}}\subset T$ and also  
$g_{r}x_{r}g_{r}^{-1}\in T$. By \cite[Lemma IV 2.5]{BD} for each 
$1\le r\le k$ we can find some $w_{r}\in W$ such that 
$w_{r}x_{r}=g_{r}x_{r}g_{r}^{-1}$. We conclude then that 
$S_{r,0}=w_{r}T_{I_{r},0}$ and this in turn implies that 
$S_{r}=w_{r}T_{I_{r}}=T_{w_{r}^{-1}I_{r}}$ for all $1\le r\le k$. This proves that 
any chain $\Si_{\i}=\{S_{0}\subset\cdots\subset S_{k}\}$  
in $\T(G)$ is conjugated to a chain of the form 
$T_{J_{0}}\subset T_{J_{1}}\subset \cdots \subset T_{J_{k}}$ 
for a collection of sets of roots $J_{0},\dots, J_{k}$. Moreover, 
two chains such chains $T_{I_{0}}\subset T_{I_{1}}\subset \cdots \subset T_{I_{k}}$  
and $T_{J_{0}}\subset T_{J_{1}}\subset \cdots \subset T_{J_{k}}$ are 
conjugated if and only if we can find some $w\in W$ such that 
$T_{I_{r}}=T_{wJ_{r}}$ for $0\le r\le k$.  This proves that the set 
$\E(\i)$ can be identified with the set of equivalence classes 
of sequences of sets of roots of the form 
$(J_{0},\dots, J_{k})$ with $\rho(T_{J_{r}})=i_{r}$ for $0\le r\le k$, where 
$(J_{0},\dots, J_{k})\sim (I_{0},\dots, I_{k})$ if and only if 
we can find some $w\in W$ such that $T_{I_{r}}=T_{wJ_{r}}$ for $0\le r\le k$.

\begin{example}\label{case SU(2)} Suppose that $G=SU(2)$ which 
is a Lie group of rank $1$. In this 
case the poset $\T(G)$ has one element corresponding to the 
center of $G$ (isomorphic to $\Z/2$) which has rank zero and 
an element for every maximal torus in $G$. 
Fix $T\subset G$ the maximal torus consisting of those 
$2\times 2$ diagonal matrices in $G$. The Weyl group 
$W=\Z/2$ acts by permutation on the diagonal 
entries for such matrices. The space of all maximal tori in $G$ is 
homeomorphic to $G/N_{G}(T)$. Therefore 
$\T(G)=*\sqcup G/N_{G}(T)\cong *\sqcup \R P^{2}$.
By Theorem \ref{replacing discrte homotopy colims} it follows that 
$B_{com}G\simeq \hocolim_{\i\in \So(1)}\Fo_{G}(\i)$. 
In this case we have 
\begin{align*}
\Fo_{G}(0)&=BZ(G)= B\Z/2= \R P^{\infty},\\
\Fo_{G}(1)&=G/T\times_{W} BT
=\SS^{2}\times_{\Z/2}\C P^{\infty},\\
\Fo_{G}(0,1)&=G/N_{G}(T)\times B\Z/2
= \R P^{2}\times \R P^{\infty}.
\end{align*}
Therefore $B_{com}G$ is homotopy equivalent to 
the homotopy pushout of the diagram
\[
\Fo_{G}(0)\stackrel{p_{0}}{\longleftarrow}  
\Fo_{G}(0,1)\stackrel{p_{1}}{\longrightarrow} \Fo_{G}(1),
\] 
where 
\[
p_{0}:\Fo_{G}(0,1)\cong
\R P^{2}\times \R P^{\infty}\to \Fo_{G}(0)\cong\R P^{\infty}
\] 
corresponds to second projection  and 
\[
p_{1}:\Fo_{G}(0,1)\cong \R P^{2}\times \R P^{\infty}\to 
\Fo_{G}(1)\cong \SS^{2}\times_{\Z/2}\C P^{\infty}
\] 
is the map induced by the inclusion $\Z/2\hookrightarrow T$. Using 
the associated Mayer--Vietoris sequence we obtain
\begin{equation*}
H^{k}(B_{com}SU(2);\Z)\cong \left\{ 
\begin{array}{ccl}
\Z &  \text{ if } & k=0, \\
0 &  \text{ if } & k=2 \text{ or $k$ odd}, \\
\Z\oplus\Z &  \text{ if } & k\equiv 0 \text{ (mod 4)}, k>0, \\
\Z/2 &  \text{ if } & k\equiv 2 \text{ (mod 4)}, k>2. \\
\end{array}
\right. 
\end{equation*}

Note in particular the presence of a  $\Z/2$-factor in 
cohomological degrees $k\equiv 2 \text{ (mod 4)}$ 
and  $k>2$.  The existence 
of this $2$-torsion is particularly 
intriguing and we are interested in 
finding a suitable geometric interpretation.
\end{example}

Suppose now that $G$ is a Lie group that is  compact 
and connected Lie group 
with center $Z$. Let  $n={\rm{rank}}(G)-{\rm{rank}}(Z)$. 
Define a functor $\H_{G}:\So(n)\to \Top$ as follows. If 
$\i=\{i_{0},\dots,i_{k}\}$ is an object in  $\So(n)$ then we define 
\[
\H_{G}(\i):=\{(S_{0},\dots, S_{k},x) 
~|~ S_{0}\subset\cdots\subset S_{k}\in \T(G), \  
\rho(S_{r})=i_{r} \text{ for } 0\le r\le k 
\text{ and } x\in G/S_{0} \}. 
\]
If $\j$ is a subset of $\i$ then the corresponding function is the map 
$p_{\i,\j}:\H_{G}(\i)\to \H_{G}(\j)$ induced by the projection maps 
and the quotient map. 

\begin{theorem}\label{homotopy colim E_{com}}
Suppose that $G$ is a compact connected Lie group. Then there is a 
natural $G$-equivariant homotopy equivalence
\[
\hocolim_{\i\in \So(n)}\H_{G}(\i)\simeq  E_{com}G_{\BONE}.
\]
\end{theorem}
\Proof
To start  define a functor $\tilde{\H}_{G}:\So(n)\to \Top$ that 
for an object $\i$ in $\So(n)$ corresponds  
\[
\tilde{\H}_{G}(\i):=\{(S_{0},\dots, S_{k},x) 
~|~ S_{0}\subset\cdots\subset S_{k}, \  
\rho(S_{r})=i_{r} \text{ for } 0\le r\le k 
\text{ and } x\in p_{com}^{-1}(BS_{0}) \}, 
\]
where $p_{com}:E_{com}G_{\BONE}\to B_{com}G$ is the projection map. 
Note that $p_{com}^{-1}(BS_{0})$ is the geometric realization of the
subsimplicial space of $[E_{com}G]_*$ whose $n$-th space is 
$\Hom(\Z^{n},S_0)\times G=S_0^{n}\times G$. In particular  
 we have a homotopy equivalence 
 $p_{com}^{-1}(BS_{0})\simeq ES_0\times_{S_0}G\simeq G/S_{0}$.
Using this equivalence we obtain a 
natural transformation $\mu:\H_{G}\to \tilde{\H}_{G}$ 
such that $\mu_{\i}:\H_{G}(\i)\to \tilde{\H}_{G}(\i)$ is a  
$G$-equivariant homotopy 
equivalence for every $\i$. We conclude then that there is a 
$G$-equivariant homotopy equivalence 
\[
\hocolim_{\i\in \So(n)}\H_{G}(\i)\simeq \hocolim_{\i\in \So(n)}\tilde{\H}_{G}(\i).
\]
On the other hand, let $\Do$ be the topological category whose objects are pairs 
of the form $(S,x)$, where $S\in \T(G)$ and $x\in p_{com}^{-1}(BS)$. There is a unique 
morphism  $(S_{1},x)\to (S_{2},y)$ in $\Do$ if and only if $x=y$ and
$S_{1}\subset S_{2}$. An argument similar to the one provided in 
Theorem \ref{replacing discrte homotopy colims} shows that there is a
$G$-equivariant homeomorphism $\hocolim_{\i\in \So(n)}\H_{G}(\i)\cong B\Do$. 
Finally,  let $\Co$ be the topological category whose objects are 
the elements in $E_{com}G_{\BONE}$ and the only morphisms are the 
identity morphisms. Thus we have   $B\Co=E_{com}G_{\BONE}$.  
Let $F:\Do\to \Co$ be the functor 
that for an object $(S,x)$ in $\Do$ corresponds 
$F(S,x)=x\in p_{com}^{-1}(BS)\subset E_{com}G_{\BONE}$. The functor 
$F$ sends any morphism in $\Do$ to the corresponding identity morphism 
in $\Co$. Using the same argument as in 
Theorem  \ref{replacing discrte homotopy colims} 
we conclude that the map $BF:B\Do\to B\Co$ is a $G$-equivariant 
homotopy equivalence.
\qed

\begin{remark}
Let $Y(G):=\hocolim_{\i\in \So(n)}\H_{G}(\i)$. Note that $Y(G)$ 
is a finite $G$--CW-complex and by Theorem \ref{homotopy colim E_{com}} 
we have $B_{com}G_{\BONE}\simeq EG\times_{G} Y(G)$.  
When $G$ is abelian $Y(G)$ is a contractible space, 
and so it can be seen as measuring how far $G$ is from being a 
commutative group.  Using the Atiyah-Segal 
completion theorem we conclude that  
$K^{*}(\BC G_{\BONE})$ is the completion of $K^{*}_{G}(Y(G))$ with 
respect to the augmentation ideal $I_{G}$ in the complex 
representation ring $R(G)$ of $G$. This 
can be seen as a generalization of the classical computation 
$K^{*}(BG)\cong R(G)^{\wedge}_{I_{G}}$ for a compact Lie group $G$. 
\end{remark}

\section{Rational cohomology of $B_{com}G_{\BONE}$}
\label{rational cohomology}

In this section we provide computations for the 
cohomology of the spaces $\BC G_{\BONE}$ 
with rational coefficients for a real or complex reductive algebraic 
group $G$ that is connected as a topological space. 
By Theorem \ref{reduction to compact} we can 
work with compact connected Lie groups without loss of 
generality. Throughout this section we take the rational numbers 
as the ground field for all computations unless otherwise 
specified.

\medskip

Fix a compact connected 
Lie group $G$. Let $T\subset G$ be a maximal torus and let $W$ 
be the corresponding Weyl group. For every $n\ge 0$ consider 
the map
\begin{align*}
\bar{\varphi_{n}}:G\times T^{n}&\to \Hom(\Z^{n},G)_{\BONE}\\
(g,t_{1},\dots,t_{n})&\mapsto (gt_{1}g^{-1},\dots,gt_{n}g^{-1}).
\end{align*}
The group $N_{G}(T)$ acts naturally 
on $G\times T^{n}$ by
\[
n\cdot (g,t_{1},\dots,t_{n}):=(gn^{-1},nt_{1}n^{-1},\dots,nt_{n}n^{-1}).
\]
The map $\bar{\varphi_{n}}$ is invariant under this action, as a 
result we obtain a continuous map 
\begin{align*}
\varphi_{n}:G/T\times_{W} T^{n}=G\times_{N_{G}(T)} T^{n}
&\to \Hom(\Z^{n},G)_{\BONE}\\
[(g,t_{1},\dots,t_{n})]&\mapsto (gt_{1}g^{-1},\dots,gt_{n}g^{-1}).
\end{align*}
Here $W$ acts diagonally on $T^{n}$. This map is surjective as 
any $n$-tuple 
$(g_{1},\dots,g_{n})$ of elements in $G$ belongs to 
$\Hom(\Z^{n},G)_{\BONE}$ if and only there is a maximal 
torus in $G$ that contains $g_{1}\dots,g_{n}$ and all maximal 
tori in $G$ are conjugated. By \cite[Lemma 3.2]{Baird} 
it follows that the fibers of $\varphi_{n}$ are rationally acyclic 
and thus $\varphi_{n}$ induces an isomorphism in cohomology 
with rational coefficients.
It is easy to see that the collection 
$\{\varphi_{n}\}_{n\ge0}$ defines a map of simplicial spaces and  
by passing to the geometric realization we obtain 
a continuous surjective map
\[
\varphi:G/T\times_{W}BT\to \BC G_{\BONE}.
\]
In the same way as in  \cite[Theorem 6.1]{ACT},  
we conclude that the map $\varphi$  induces an isomorphism in 
cohomology with rational coefficients and thus we 
obtain an isomorphism
\begin{equation}\label{theorem ACT}
\varphi^{*}:H^{*}(\BC G_{\BONE})\stackrel{\cong}
{\rightarrow}\left(H^{*}(G/T)\otimes H^{*}(BT) \right)^{W},
\end{equation}
with $W$ acting diagonally on $H^{*}(G/T)\otimes H^{*}(BT) $.
This can be used to provide  the following useful identification of 
the rational cohomology of $\BC G_{\BONE}$.

\begin{proposition}\label{interchanging}
Suppose that $G$ is a compact connected Lie group and 
let $T\subset G$ be a maximal torus. 
Then there is a natural isomorphism of  rings 
\[
\alpha_{G}: H^{*}(\BC G_{\BONE})\stackrel{\cong}
{\rightarrow} \left (H^{*}(BT)\otimes 
H^{*}(BT)\right)^{W}/J_{G}.
\]
In the above equation  
$W$ acts diagonally on $H^{*}(BT)\otimes H^{*}(BT)$ 
and $J_{G}$ is the ideal generated by the 
elements of positive degrees in the image of
\begin{align*}
i_{1}:H^{*}(BG)&\to \left(H^{*}(BT)\otimes H^{*}(BT)\right)^{W}\\
x&\mapsto x\otimes 1.
\end{align*}
\end{proposition}

\Proof
 The Eilenberg-Moore spectral sequence with 
$\Q$-coefficients associated to the fibration  
$G/T\to BT\stackrel{i}{\rightarrow} BG$ collapses at the 
$E_{2}$-term (see \cite{mccleary}, page 278).
Also, there is a $W$-equivariant isomorphism of graded 
rings  $H^{*}(G/T)\cong H^{*}(BT)/I_{G}$,
where $I_{G}$ is the ideal in $ H^{*}(BT)$ generated 
by the elements of positive degree in the image of 
$i^{*}:H^{*}(BG)\to H^{*}(BT)$. 
Consider now the natural map 
\[
\pi:H^{*}(BT)\otimes H^{*}(BT)\to 
(H^{*}(BT)/I_{G})\otimes H^{*}(BT).
\]
This is a surjective map whose kernel is the ideal $\tilde{I}_{G}$ 
generated by the elements of positive degree in the image of the 
map $i_{1}:H^{*}(BG)\to H^{*}(BT)\otimes H^{*}(BT)$ given by
$x\mapsto x\otimes 1$.
Thus we have a short exact sequence 
\[
0\to\tilde{I}_{G} \to  H^{*}(BT)\otimes H^{*}(BT) 
\to (H^{*}(BT)/I_{G})\otimes H^{*}(BT)\to 0.
\]
Since we are working in characteristic zero and $W$ is a 
finite group,  the exactness 
of this sequence is preserved at the level of $W$-invariants; that is, 
there is a short exact sequence
\begin{equation}\label{W-invariants}
0\to \tilde{I}_{G}^{W}\to \left( H^{*}(BT)\otimes H^{*}(BT)\right)^{W}
\to \left(  (H^{*}(BT)/I_{G})\otimes H^{*}(BT)\right)^{W}\to 0.     
\end{equation}
Note that $J_{G}=\tilde{I}_{G}^{W}$, thus we obtain a natural 
isomorphism 
\[
\psi:\left( H^{*}(BT)\otimes H^{*}(BT)\right)^{W}/J_{G}\to 
\left(H^{*}(G/T)\otimes H^{*}(BT) \right)^{W}.
\]
The required isomorphism is then $\alpha_{G}:=\psi^{-1}\circ \varphi$, 
where $\varphi$ is as in (\ref{theorem ACT}).
\qed

\medskip

The previous proposition has a number of interesting applications. 
To start note that we have a natural inclusion 
$BT\subset \BC G_{\BONE}\subset BG$.  At the level of cohomology 
groups this induces a natural monomorphism 
$H^{*}(BG)\hookrightarrow H^{*}(\BC G_{\BONE})$ and thus we can consider
$H^{*}(\BC G_{\BONE})$ as a module over $H^{*}(BG)$. 
Under the identification 
\[
\alpha: H^{*}(\BC G_{\BONE})\stackrel{\cong}
{\rightarrow} \left (H^{*}(BT)\otimes 
H^{*}(BT)\right)^{W}/J_{G}
\]
provided in the previous proposition, this structure as 
$H^{*}(BG)$--module corresponds to the structure on 
$\left (H^{*}(BT)\otimes H^{*}(BT)\right)^{W}/J_{G}$ given by 
$g\cdot [x\otimes y ]:=[x\otimes gy]$. As a consequence of this 
we derive the following theorem.

\begin{theorem}
Suppose that  $G$ is a compact, connected Lie group.
Then 
$H^{*}(\BC G_{\BONE})$ is a free module over $H^{*}(BG)$ 
of rank $|W|$, where $W$ is the corresponding Weyl group.
\end{theorem}

\Proof 
Fix a maximal torus $T\subset G$  and let $W$ be the 
corresponding Weyl group. Let $S=H^{*}(BG)$ 
and $A=H^{*}(BT)$. These are graded rings, $W$ 
acts on $A$ with degree-preserving ring automorphisms
and we  have a natural isomorphism $S\cong A^{W}$. The ring 
$S$  is a polynomial ring 
on finitely many commuting variables. Also, the ring
$A$ can be seen as a graded module over $S$ and this is 
in fact a  free module  of rank $|W|$. 
Consider $M^{W}:=(A\otimes A)^{W}$. This 
is a graded ring that contains $R:=S\otimes S$ as a subring. Thus 
$M^{W}$ can be seen as a graded module over $R$.
As a first step we will show that $M^{W}$ is a finitely generated free 
$R$-module.  To this end, note that $R$ 
is a Cohen-Macaulay ring as it is a polynomial ring over $\Q$. 
The same is true for $A\otimes A$. 
Since $W$ acts by degree--preserving 
ring automorphisms on  $A\otimes A$, it follows that 
$M^{W}=(A\otimes A)^{W}$ is also a Cohen-Macaulay ring by 
\cite[Proposition 13]{HE}. We observe that $M^{W}$ is finitely 
generated as an $R$-module. 
Indeed, suppose  $\{e_{w}\}_{w\in W}$ 
is a free basis of $A$ as a module over $S$. 
Then $\{e_{v}\otimes e_{w}\}_{v,w\in W}$ is a free basis of  
$A\otimes A$ as a module over $R=S\otimes S$. 
Define the averaging operator
\begin{align*}
\rho:A\otimes A&\to  (A\otimes A)^{W}=M^{W}\\
f&\mapsto \frac{1}{|W|}\sum_{w\in W}w\cdot f.
\end{align*}
The map $\rho$ is surjective and $R$-linear. 
This implies that the collection 
$\{\rho(e_{v}\otimes e_{w})\}_{v,w\in W}$ generates 
$M^{W}$ as a module over $R$.  Thus $M^{W}$ is a 
finitely generated $R$-module. Since $R$ is 
a polynomial algebra over $\Q$, then any finitely generated 
$R$-module has finite projective dimension.
Using the Auslander-Buchsbaum formula for graded rings we 
obtain 
\[
{\rm{pd}}_{R}(M^{W})={\rm{depth}}(R)
-{\rm{depth}}(M^{W}),
\]
where ${\rm{pd}}_{R}(M^{W})$ is the projective 
dimension of $M^{W}$ as an $R$-module.
Since both $M^{W}$ and $R$ are Cohen-Macaulay rings 
and $M^{W}$ is finitely generated 
as an $R$-module,  this implies  
${\rm{pd}}_{R}(M^{W})={\rm{dim}}(R)
-{\rm{dim}}(M^{W})=0$. 
This means that  $M^{W}$ is projective as an 
$R$-module and by the Quillen-Suslin theorem (see 
\cite[Theorem 4]{Quillen1}),
it follows that $M^{W}$ is a free $R$-module. 
Recall that $J_{G}$ is the ideal in $M^{W}$
generated by the elements of positive degree of the 
form $x\otimes 1$ for $x\in S$. 
Suppose that $\{a_{w}\}_{w}$ is a free basis of $M^{W}$ 
as a module over $R$. If $f_{w}=\bar{a}_{w}$ is the image 
of $a_{w}$ in $M^{W}/J_{G}$ under the natural map, then 
it follows that $\{f_{w}\}_{w}$ is a free basis of 
$M^{W}/J_{G}$ as a module 
over $S$.  By Theorem \ref{interchanging}
this means that $H^{*}(\BC G_{\BONE})$ is  free 
as a module over $H^{*}(BG)$.  To finish we only need to 
compute the rank of $H^{*}(\BC G_{\BONE}) $. For this 
recall that we have a natural isomorphism 
$\varphi^{*}:H^{*}(\BC G_{\BONE})\stackrel{\cong}
{\rightarrow}\left(H^{*}(G/T)\otimes H^{*}(BT) \right)^{W}$.
As an ungraded $W$-module $H^{*}(G/T)$ is isomorphic to the 
regular $W$-representation. It follows that 
 as an ungraded module, 
$H^{*}(\BC G_{\BONE})$ is isomorphic to $H^{*}(BT)$ and 
the latter has rank $|W|$ as a module over $H^{*}(BG)$. This 
implies that as a graded  $H^{*}(BG)$-module 
$H^{*}(\BC G_{\BONE})$ is free and of rank $|W|$. 
 \qed

\begin{remark} The previous theorem is not true 
in general if we use integer coefficients. For example  if
$G=SU(2)$ then $H^{*}(B_{com}SU(2);\Z)$ is not free as a 
module $H^{*}(BSU(2);\Z)$ 
because the former contains $2$-torsion as we proved in 
Example \ref{case SU(2)} and 
the latter does not contain torsion.
\end{remark}

Consider now the inclusion map $i:\BC G_{\BONE}\to BG$. 
Up to homotopy we have a fibration sequence 
\begin{equation}\label{fibration \BC G}
E_{com}G_{\BONE}\stackrel{p_{com}}
{\rightarrow} \BC G_{\BONE}\stackrel{i}
{\rightarrow} BG.
\end{equation}
Since $G$ is assumed to be connected then 
the base space, $BG$, 
is simply connected. The $E_{2}$-term of the Eilenberg-Moore 
spectral sequence with $\Q$-coefficients 
associated to the fibration (\ref{fibration \BC G}) 
is 
\[
E_{2}^{*,*}=\Tor_{H^{*}(BG)}(\Q,H^{*}(\BC G_{\BONE}))
\]
and this spectral sequence converges strongly to 
$H^{*}(E_{com}G_{\BONE})$. By the previous theorem, 
if $G$ is a compact connected Lie group 
then $H^{*}(\BC G_{\BONE})$ is a 
free module over $H^{*}(BG)$. It follows that 
\[
\Tor_{H^{*}(BG)}(\Q,H^{*}(\BC G_{\BONE})) 
\cong\Q\otimes_{H^{*}(BG)}H^{*}(\BC G_{\BONE})
\]
and the Eilenberg-Moore spectral sequence collapses to the 
$E_{2}^{0,*}$-column. The map 
$p_{com}^{*}:H^{*}(\BC G_\BONE{})\to H^{*}(E_{com}G_{\BONE})$
is surjective since $\text{Im}(p_{com}^*)=E_{\infty}^{0,*}$ and the 
sequence collapses at the $E_{2}^{0,*}$-column. 
Let $K_{G}$ denote the ideal in $H^{*}(\BC G_{\BONE})$ 
generated by the elements in $H^{*}(BG)$ of positive degree. 
Then $\text{Ker}(p_{com}^*)=K_{G}$ and we conclude 
that there is an isomorphism of rings 
$H^{*}(E_{com}G_{\BONE})\cong H^{*}(\BC G_{\BONE})/K_{G}$.
Using the isomorphism provided in Theorem \ref{interchanging} 
we conclude that if
$L_{G}$ is the ideal in
$(H^{*}(BT)\otimes H^{*}(BT))^{W}$ generated by the
elements of positive degree in the image of
$H^{*}(BG)\otimes H^{*}(BG)$, then 
there is a natural isomorphism 
$H^{*}(E_{com}G_{\BONE})\cong 
(H^{*}(BT)\otimes H^{*}(BT))^{W}/L_{G}
\cong (H^*(G/T)\otimes H^*(G/T))^W$.
This proves the following corollary. 

\begin{corollary}\label{Poincare series}
Suppose that $G$ is a connected compact Lie group with
maximal torus $T$ and associated Weyl group $W$.
Then there is a natural 
isomorphism of rings 
\[
H^{*}(E_{com}G_{\BONE})
\cong (H^*(G/T)\otimes H^*(G/T))^W,
\]
and the Poincar\'e series of $\BC G_{\BONE}$ and 
$E_{com}G_{\BONE}$ satisfy 
$P_{\BC G_{\BONE}}(t)=P_{BG}(t)P_{E_{com}G_{\BONE}}(t)$.
\end{corollary}

Note that $G/T\times G/T$ is a compact, orientable manifold and that
$W$ preserves orientation. Hence we infer that the fundamental class
is $W$--invariant. This yields the following:

\begin{corollary}
The following three statements are equivalent for a compact
connected Lie group $G$:
(1)~$E_{com}G_{\BONE}$ is contractible;
(2)~$E_{com}G_{\BONE}$ is rationally acyclic; and
(3)~$G$ is abelian.
\end{corollary}

\Proof
If $E_{com}G_{\BONE}$ is contractible then it's acyclic. If
it's acyclic, then $G/T$ must be zero--dimensional, hence $G=T$
and so $G$ is abelian. If $G$ is abelian $\BC G_{\BONE}=BG$ and so
$E_{com}G_{\BONE}$ is contractible.
\qed

\medskip

From the description given in Corollary \ref{Poincare series},  
it follows that the Poincar\'e series for $E_{com}G_{\BONE}$ 
encodes information about all  the
complex irreducible  representations of the Weyl group 
$W$ associated to the pair $(G,T)$. To see this   
recall that as an ungraded $W$--representation 
$H^{*}(G/T;\C)$ is isomorphic to the regular representation 
and thus it contains all the irreducible representations of $W$.
For every irreducible representation $\lambda$ of $W$, 
consider its character $\chi^{\lambda}$. 
The multiplicity of $\lambda$ in the regular 
representation equals its degree which we denote by 
$f^{\lambda}=\chi^{\lambda}(e)$.  
The multiplicity of the representation 
$\lambda$ in the representations  $H^{i}(G/T;\C)$ for 
$i\ge 0$ can be described by the fake degree polynomial 
$f^{\lambda}(t)$ defined by 
\[
f^{\lambda}(t):=
\sum_{i\ge 0 } t^{i} \left<\chi^{\lambda},H^{i}(G/T;\C)  \right>.
\]
In other words, the coefficient of $t^{i}$ in $f^{\lambda}(t)$ 
is exactly the multiplicity of $\lambda$ in $H^{i}(G/T;\C)$.
The Poincar\'e polynomial of the flag manifold $G/T$ is then 
given by 
\[
P_{G/T}(t)=\sum_{\lambda}f^{\lambda}f^{\lambda}(t)
=\sum_{\lambda}f^{\lambda}(1)f^{\lambda}(t),
\]
where $\lambda$ runs through all complex irreducible 
representations of $W$.  On the other hand, 
the Poincar\'e polynomial of $E_{com}G_{\BONE}$ 
is given by 
\[
P_{E_{com}G_{\BONE}}(t)=
\sum_{\lambda}f^{\bar{\lambda}}(t)f^{\lambda}(t),
\]
where $\lambda$ runs through all complex irreducible 
representations of $W$, and ${\bar{\lambda}}$ is the complex 
conjugate of $\lambda$. To see this note that by 
Corollary \ref{Poincare series} and the universal coefficient
theorem we have $H^{*}(E_{com}G_{\BONE};\C)\cong 
(H^*(G/T;\C)\otimes H^*(G/T;\C))^W$. On the other hand, 
for each $0\le k\le n$ 
we have an isomorphism
\[
(H^{k}(G/T;\C)\otimes H^{n-k}(G/T;\C))^W\cong \Hom_{W}(
\Hom(H^{k}(G/T;\C),\C),H^{n-k}(G/T;\C)).
\]
This together with Schur's lemma shows that 
$H^{n}(E_{com}G_{\BONE};\C)$ is a vector space over $\C$ of 
dimension 
\[
\sum_{0\le k\le n}\sum_{\lambda} \left<\chi^{\bar{\lambda}},H^{k}(G/T;\C) \right>
 \left<\chi^{\lambda},H^{n-k}(G/T;\C) \right> 
\]
and thus $P_{E_{com}G_{\BONE}}(t)=
\sum_{\lambda}f^{\bar{\lambda}}(t)f^{\lambda}(t)$.

\section{The cases $SU(n)$, $U(n)$ and 
$Sp(n)$ }

In this section we study in detail the cohomology with rational 
coefficients of the space $B_{com}G$ when $G$ is one of the 
classical groups $SU(n)$, $U(n)$ and 
$Sp(n)$ and also for their corresponding complexifications 
$SL_{n}(\C)$, $GL_{n}(\C)$ and  $Sp_{n}(\C)$.  In particular 
we provide explicit free bases of $H^{*}(B_{com}G;\Q)$ as 
a module over $H^{*}(BG;\Q)$. As in the previous section, 
we take the rational numbers as the ground field unless 
otherwise specified. 

\medskip

To start recall that by \cite[Proposition 2.5]{AG}  
we have that $\Hom(\Z^{n},G)$ is path--connected for all 
$n\ge 0$ when $G$ is one of the groups $SU(n)$, $U(n)$ and 
$Sp(n)$. Thus for such groups $G$ and their corresponding
complexifications we have 
$\BC G=\BC G_{\BONE}$ and $E_{com}G=E_{com}G_{\BONE}$. 

\subsection{Case $G=U(n)$} Suppose $G=U(n)$. 
In this case we can choose 
a maximal torus $T\subset U(n)$ to be the set 
of diagonal matrices with entries in $\SS^{1}$. We have 
$H^{*}(BT)\cong \Q[\x]$, where $\x:=\{x_{1},\dots, x_{n}\}$ and 
$\text{deg}(x_{i})=2$ for $1\le i\le n$. The Weyl group 
is the symmetric group 
$W=\Sigma_{n}$ acting by permutation on the  variables 
$x_{1},\dots, x_{n}$.  Therefore by Theorem \ref{interchanging} 
we have an isomorphism 
\[
\alpha_{n}:=\alpha_{U(n)}:H^{*}(B_{com}U(n))\stackrel{\cong}
{\rightarrow}
(\Q[\x]\otimes \Q[\y])^{\Sigma_{n}}/J_{U(n)},
\]
where $\Sigma_{n}$ acts diagonally by permuting the variables 
$\x=\{x_{1},\dots,x_{n}\}$ and $\y=\{y_{1},\dots,y_{n}\}$. The algebra
$M^{\Sigma_{n}}:=(\Q[\x]\otimes \Q[\y])^{\Sigma_{n}}$
is known as the algebra of multisymmetric polynomials. In this 
case $J_{n}:=J_{U(n)}$ is the ideal in $M^{\Sigma_{n}}$ 
generated by the elementary 
symmetric polynomials 
\[
e_{k}(x_{1},\dots,x_{n})=\sum_{1\le i_{1}<i_{2}<\cdots <i_{k}\le n}
x_{i_{1}}x_{i_{2}}\cdots x_{i_{k}}
\]
for $1\le k\le n$.  Since we are working with rational coefficients, 
the ideal $J_{n}$ is also the ideal generated by the power sums
$p_{n}(a,0):=x_{1}^{a}+\cdots +x_{n}^{a}$
for $1\le a\le n$. These classical power sums have analogues in the 
ring of multisymmetric polynomials. For every pair of integers 
$a,b\ge 0$ define the power sum 
$p_{n}(a,b):=x_{1}^{a}y_{1}^{b}+
\cdots+x_{n}^{a}y_{n}^{b}$.
Clearly $p_{n}(a,b)\in M^{\Sigma_{n}}$ for all $a,b\ge 0$. Moreover, 
it is well known that the polynomials $p_{n}(a,b)$ for $0<a+b\le n$ 
generate the $\Q$-algebra $M^{\Sigma_{n}}$ although these 
polynomials are not algebraically independent.  
(See for example \cite{Vaccarino} for a modern 
account on multisymmetric polynomials).  We know by 
Theorem \ref{free module} that $M^{\Sigma_{n}}/J_{n}$ is a free 
module over $H^{*}(BU(n))\cong \Q[\y]^{\Sigma_{n}}$. 
An explicit free basis for 
$M^{\Sigma_{n}}/J_{n}$ as a module over 
$\Q[\y]^{\Sigma_{n}}$ can be constructed using the work 
in \cite{Allen}. For this consider the averaging operator 
\begin{align*}
\rho:\Q[\x]\otimes \Q[\y]
&\to  (\Q[\x]\otimes \Q[\y])^{\Sigma_{n}}=M^{\Sigma_{n}}\\
f(\x,\y)&
\mapsto \frac{1}{n!}\sum_{w\in\Sigma_{n}} 
f(w\x,w\y).
\end{align*} 
For every $w\in \Sigma_{n}$ the diagonal descent monomial 
is defined to be
\[
e_{w}:=\prod_{w^{-1}(i)>w^{-1}(i+1)}(x_{1}\cdots x_{i}) \otimes 
\prod_{w(j)>w(j+1)}(y_{w(1)}\cdots y_{w(j)}).
\]
By \cite[Theorem 1.3]{Allen} the collection 
$\{\rho(e_{w})\}_{w\in \Sigma_{n}}$ forms a free basis of 
$M^{\Sigma_{n}}$ as a module over 
$\Q[\x]^{\Sigma_{n}}\otimes  \Q[\y]^{\Sigma_{n}}$. 

\begin{example}
Suppose that $n=3$. In this case we obtain the basis 
following basis of  $M^{\Sigma_{3}}$ as a module over 
$\Q[\x]^{\Sigma_{3}}\otimes  \Q[\y]^{\Sigma_{3}}$
\[
e_{1}=1, e_{2}=\rho(x_{1}y_{2}), e_{3}=\rho(x_{1}y_{2}y_{3}), 
e_{4}=\rho(x_{1}x_{2}y_{3}), e_{5}=\rho(x_{1}x_{2}y_{1}y_{3}), 
e_{6}=\rho(x_{1}^{2}x_{2}y_{3}^{2}y_{2}).
\]
\end{example}

Let $f_{w}$ be the 
image of $\rho(e_{w})$ in $M^{\Sigma_{n}}/J_{n}$. Then 
it follows that
$\{f_{w}\}_{w\in \Sigma_{n}}$ forms a free basis of 
$H^{*}(B_{com}U(n))\cong M^{\Sigma_{n}}/J_{n}$ as a module over 
$H^{*}(BU(n))\cong \Q[\y]^{\Sigma_{n}}$. 
For each $w\in \Sigma_{n}$ define  the descent of $w$ to be 
the set
\[
Des(w):=\{1\le i\le n-1 ~|~ w(i)>w(i+1) \}.
\]
The major index of $w$, denoted by $\maj(w)$, is defined to be 
\[
\maj(w):=\sum_{i\in Des(w)}i=\sum_{w(i)>w(i+1)}i. 
\]
For every $w\in \Sigma_{n}$ we have 
${\rm{deg}}f_{w}=2(\maj(w)+\maj(w^{-1}))$. As a 
corollary the following is obtained.

\begin{corollary}
Suppose that $n\ge 1$.  Then
\begin{align*}
P_{E_{com}GL_{n}(\C)}&=P_{E_{com}U(n)}(t)
= \sum_{w\in \Sigma_{n}}t^{2(\maj(w)+\maj(w^{-1}))},\\
P_{B_{com}GL_{n}(\C)}&=P_{B_{com}U(n)}(t)
=\frac{\sum_{w\in \Sigma_{n}}
t^{2(\maj(w)+\maj(w^{-1}))}}{\prod_{1\le i \le n}(1-t^{2i})}.
\end{align*}
\end{corollary}

Consider now he standard inclusion $i_{n}:U(n)\to U(n+1)$.
Let $U:=\colim_{n\to \infty}U(n)$. Then 
$B_{com}U=\colim_{n\to \infty}B_{com}U(n)$ 
and $H^{*}(B_{com}U(n))\cong \varprojlim H^{*}(B_{com}U(n))$.
The isomorphisms $\alpha_{n}$ and the maps $i_{n}$ 
are compatible in the sense that 
\begin{equation*}
\begin{CD}
H^{*}(B_{com}U(n+1))@>{\alpha_{n+1}}>>M^{\Sigma_{n+1}}/J_{n+1}\\
@V{i_{n}^{*}}VV     @VV{j_{n}^{*}}V\\
H^{*}(B_{com}U(n))@>{\alpha_{n}}>>
M^{\Sigma_{n}}/J_{n}.\\
\end{CD}
\end{equation*}
is a commutative diagram, where 
$j^{*}_{n}:M^{\Sigma_{n+1}}\to M^{\Sigma_{n}}$ is the map obtained by 
sending $x_{i}\mapsto x_{i}$, $y_{i}\mapsto y_{i}$  for $1\le i\le n$ and 
$x_{n+1}, y_{n+1}\mapsto 0$. Define
\[
M^{\Sigma_{\infty}}:= \varprojlim M^{\Sigma_{n}} \text{ and }
J_{\infty}:= \varprojlim J_{n}.
\]
Then we obtain an isomorphism of graded $\Q$-algebras
\[
H^{*}(B_{com}U)\cong \varprojlim M^{\Sigma_{n}} J_{n}
\cong M^{\Sigma_{\infty}}/J_{\infty}.
\]
The last isomorphism follows from the fact that 
$j^{*}_{n}:J_{n+1}\to J_{n}$ is surjective for every $n\ge 1$ 
and thus $ \varprojlim^{1} J_{n}$ vanishes. 
We show next that this algebra is a polynomial algebra.
For every pair of integers $a,b\ge 0$ we have 
 $j^{*}_{n}(p_{n+1}(a,b))=p_{n}(a,b)$ and thus these polynomials 
define an element in $M^{\Sigma_{\infty}}$ which we denote by $p(a,b)$. 
By \cite[Theorem 2]{Vaccarino} 
the algebra $M^{\Sigma_{\infty}}$ is a polynomial algebra over 
$\Q$ generated by the elements $p(a,b)$ for $(a,b)\ne 0$. 
On the other hand, since $J_{n}$ is the ideal 
in $M^{\Sigma_{n}}$ generated by the power sums 
$p_{n}(1,0),\dots, p_{n}(n,0)$, it follows that $J_{\infty}$ is the ideal 
generated by $p(a,0)$ for $a\ge 1$.  For each pair of integers 
$a,b\ge 0$ not both zero let $z_{a,b}$ be a $2(a+b)$-dimensional 
variable such that the collection $\{z_{a,b}\}_{a,b}$ is a collection 
of commuting independent variables. Define 
$\Ma_{U}=\{(a,b)\in \N^{2}~|~  b>0\}$;
we conclude that the assignment 
\begin{align*}
 \Q[z_{a,b}~|~ (a,b)\in\Ma_{U}]&\to M^{\Sigma_{\infty}}/J_{n}
 \cong H^{*}(B_{com}U)\\
 z_{a,b}&\mapsto p(a,b)
\end{align*}
is an isomorphism of algebras over $\Q$. Also define 
$GL_{\infty}(\C):=\colim_{n\to \infty} GL_{n}(\C)$. Since 
$B_{com}U(n)\simeq B_{com}GL_{n}(\C)$ for every $n\ge 0$, 
it follows that $B_{com}U\simeq B_{com}GL_{\infty}(\C)$.
This proves the following 
corollary.
\begin{corollary}\label{comb 1}
Let $\Ma_{U}=\{(a,b)\in \N^{2}~|~ b>0\}$. Then 
we have  isomorphisms of $\Q$-algebras 
 \[
 H^{*}(B_{com}GL_{\infty}(\C))\cong
 H^{*}(B_{com}U)\cong  \Q[z_{a,b}~|~ (a,b)\in\Ma_{U}].
 \] 
\end{corollary}

\subsection{Case $G=SU(n)$} The case of the special unitary 
groups $G=SU(n)$ can be handled 
in a similar way. In this case we can choose $T\subset SU(n)$ 
to be the set of diagonal matrices with entries in $\SS^{1}$ 
and determinant one. The Weyl group is the symmetric group 
$W=\Sigma_{n}$ acting by permutation on the diagonal entries 
and $H^{*}(BT)\cong \Q[\x]/(p_{n}(1,0))$, where as before 
we use the notation $\x=\{x_{1},\dots,x_{n}\}$.
Using an argument similar to that in Theorem \ref{interchanging}, 
we conclude that 
\[
H^{*}(B_{com}SU(n))\cong 
(\Q[\x]\otimes\Q[\y])^{\Sigma_{n}}/K_{n}
=M^{\Sigma_{n}}/K_{n},
\]
where $K_{n}$ is the ideal in $M^{\Sigma_{n}}$ generated by the 
multisymmetric polynomials  $p_{n}(a,0)$ for $1\le a\le n$ 
and $p_{n}(0,1)$. As it was pointed out before the collection 
$\{\rho(e_{w})\}_{w\in \Sigma_{n}}$ forms a free basis for 
$M^{\Sigma_{n}}$ as a module over 
$\Q[\x]^{\Sigma_{n}}\otimes \Q[\y]^{\Sigma_{n}}$. 
Let $g_{w}$ denote the image of $\rho(e_{w})$ in 
$M^{\Sigma_{n}}/K_{n}$. Then it  
follows that $\{g_{w}\}_{w\in \Sigma_{n}}$ forms a free basis for 
$H^{*}(B_{com}SU(n))\cong M^{\Sigma_{n}}/K_{n}$ as a module over 
$H^{*}(BSU(n))=
\Q[p_{n}(0,2),\dots, p_{n}(0,n)]$. 
As a corollary the following is obtained.

\begin{corollary}
Suppose that $n\ge 1$.  Then
\begin{align*}
P_{E_{com}SL_{n}(\C)}(t)&=P_{E_{com}SU(n)}(t)= 
\sum_{w\in \Sigma_{n}}t^{2(\maj(w)+\maj(w^{-1}))},\\
P_{B_{com}SL_{n}(\C)}(t)&=P_{B_{com}SU(n)}(t)=
\frac{\sum_{w\in \Sigma_{n}}t^{2(\maj(w)+\maj(w^{-1}))}}
{\prod_{2\le i \le n}(1-t^{2i})}, 
\end{align*}
where $\maj(w)$ is the major index of $w$ as defined above.
\end{corollary}

As in the case of the unitary groups we have a stabilization process 
given by the standard inclusions $i_{n}:SU(n)\to SU(n+1)$. 
Let $SU:=\colim_{n\to \infty} SU(n)$. It follows that 
\[
H^{*}(B_{com}SU)=\varprojlim H^{*}(B_{com}SU(n))
= \varprojlim  M^{\Sigma_{n}}/K_{n}.
\]
Let $K_{\infty}\subset M^{\Sigma_{\infty}}$ denote the ideal 
corresponding to the ideals  $K_{n}\subset M^{\Sigma_{n}}$ 
for $n\ge 0$. For $SU$ we have
$H^{*}(B_{com}SU)\cong M^{\Sigma_{\infty}}/K_{\infty}$.
Note that $K_{\infty}$ is precisely the ideal in $M^{\Sigma_{\infty}}$ 
generated by $p(a,0)$ for $a\ge 1$ and $p(0,1)$. Similarly define 
$SL_{\infty}(\C):=\colim_{n\to \infty} SL_{n}(\C)$. As a corollary 
the following is obtained.

\begin{corollary}\label{comb 2}
Let $\Ma_{SU}=\{(a,b)\in \N^{2}~|~ (a,b)\ne (0,1), b>0 \}$. Then 
we have isomorphisms of $\Q$-algebras 
 \[
 H^{*}(B_{com}SL_{\infty}(\C))\cong H^{*}(B_{com}SU)
 \cong \Q[z_{a,b} ~|~  (a,b)\in \Ma_{SU}].
 \] 
\end{corollary}

\subsection{Case $G=Sp(n)$} Finally, suppose that $G=Sp(n)$. 
In this case $H^{*}(BT)\cong \Q[\x]$, where 
$\x=\{x_{1},\dots,x_{n}\}$ and 
${\rm{deg}}(x_{i})=2$. The Weyl group is a semi-direct product  
$W=\Sigma_{n}\ltimes (\Z/2)^{n}$. This group can be identified with 
the group of signed permutations. More precisely, let 
\[
\Ib_{n}:=\{-n,-n+1,\dots,-1,1, \dots, n-1,n\}.
\]
Let $B_{n}$ denote the group of bijections $\sigma:\Ib_{n}\to \Ib_{n}$ 
such that $\sigma(-k)=-\sigma(k)$ for all $k\in \Ib_{n}$, with the 
composition of functions as the group operation. Under this 
identification the group $W\cong B_{n}$ acts on 
$\Q[\x]$ by signed permutations. In this case $H^{*}(BG)$  
is generated polynomial algebra on the generators
$e_{1}(x_{1}^{2},\dots, x_{n}^{2}),\dots, e_{n}(x_{1}^{2},\dots,x_{n}^{2})$, 
or equivalently, by the power sums $p_{n}(2,0),\dots, p_{n}(2n,0)$. 
Also $M^{B_{n}}:=(\Q[\x]\otimes \Q[\y])^{B_{n}}$
is the ring of diagonally signed-symmetric or signed-invariant 
multisymmetric polynomials. Note in 
particular that $M^{B_{n}}$ is a subalgebra of the algebra of 
multisymmetric polynomials 
$M^{\Sigma_{n}}$. By Theorem \ref{interchanging} 
we have an isomorphism 
\[
\alpha_{Sp(n)}:H^{*}(B_{com}Sp(n))\stackrel{\cong}
{\rightarrow} M^{B_{n}}/L_{n},
\]
where $L_{n}=J_{Sp(n)}$ is the ideal in $M^{B_{n}}$  generated 
by the power sums $p_{n}(2,0),\dots, p_{n}(2n,0)$.  An explicit 
basis for $M^{B_{n}}/L_{n}$ as a module over 
$\Q[\y]^{B_{n}}$ can be found using the work in 
\cite{Gomez}.  As before given $w\in B_{n}$ define its
descent to be the set
\[
Des(w):=\{1\le i\le n-1 ~|~ w(i)>w(i+1)\}.
\]
For $1\le i \le n $ let
\[
d_{i}(w):=|\{j\in Des(w) ~|~ j\ge i \}|,
\]
\begin{equation*}
\varepsilon_{i}(w):= \left\{
\begin{array}{ccc}
0& \text{ if } &w(i)>0,\\
1&\text{ if } &w(i)<0,
\end{array}
\right. 
\end{equation*}
and 
\[
f_{i}(w):=2d_{i}(w)+\varepsilon_{i}(w).
\]
The diagonal signed descent monomial associated 
to $w$ is defined to be
\[
c_{w}:=\left(\prod_{i=1}^{n}x_{i}^{f_{i}(w^{-1})} \right)
\left(\prod_{i=1}^{n}y_{|w(i)|}^{f_{i}(w)}
 \right)=\prod_{i=1}^{n}x_{i}^{f_{i}(w^{-1})}
 y_{i}^{f_{|w^{-1}(i)|}(w)}.
\]
By \cite[Theorem 1.1]{Gomez} the collection 
$\{\rho(c_{w})\}_{w\in \Sigma_{n}}$ forms a free basis of 
$M^{B_{n}}$ as a module over 
$\Q[\x]^{B_{n}}\otimes  \Q[\y]^{B_{n}}$.  

\begin{example}
Suppose that $n=2$. In this case we obtain the  
following basis of  $M^{B_{2}}$ as a module over 
$\Q[\x]^{B_{2}}\otimes  \Q[\y]^{B_{2}}$
\begin{align*}
c_{1}&=1, c_{2}=\rho(x_{1}y_{1}), c_{3}=\rho(x_{1}^{2}y_{2}^{2}), 
c_{4}=\rho(x_{1}y_{1}y_{2}^{2}), c_{5}=\rho(x_{1}^{2}x_{2}y_{2}),\\
c_{6}&=\rho(x_{1}x_{2}y_{1}y_{2}), c_{7}
=\rho(x_{1}^{2}x_{2}y_{1}^{2}y_{2}),
c_{8}=\rho(x_{1}^{3}x_{2}y_{1}^{3}y_{2}).
\end{align*}
\end{example}

Let  $h_{w}$ be the 
image of $\rho(c_{w})$ in $M^{B_{n}}/L_{n}$. It follows 
that  $\{h_{w}\}_{w\in B_{n}}$ forms a free basis of 
$H^{*}(B_{com}Sp(n))\cong M^{B_{n}}/L_{n}$ as a module over 
$H^{*}(BSp(n))\cong \Q[\y]^{B_{n}}$.  The flag major index of a signed permutation 
$w$ was defined in  \cite{AR} to be  
\[
\fmaj(w)=\sum_{i=1}^{n}f_{i}(w)=2\maj(w)+\text{neg}(w),
\]
where  $\maj(w):=\sum_{i\in Des(w)}i=\sum_{w(i)>w(i+1)}i$ and 
$\text{neg}(w):=|\{1\le i\le n ~|~ w(i)<0\}|$. 
Note that for every $w\in B_{n}$ we have 
${\rm{deg}}h_{w}=2(\fmaj(w^{-1})+\fmaj(w))$.
As a corollary the following is obtained.

\begin{corollary}
Suppose that $n\ge 1$.  Then
\begin{align*}
P_{E_{com}Sp_{n}(\C)}(t)&= P_{E_{com}Sp(n)}(t)
= \sum_{w\in B_{n}}t^{2(\fmaj(w^{-1})+\fmaj(w))}, \\
P_{B_{com}Sp_{n}(\C)}(t)&=P_{B_{com}Sp(n)}(t)
=\frac{\sum_{w\in B_{n}}t^{2(\fmaj(w^{-1})+\fmaj(w))}}
{\prod_{1\le i \le n}(1-t^{4i})}.
\end{align*}
\end{corollary}

As in the case of the unitary groups we have a stabilization process 
given by the standard inclusions $i_{n}:Sp(n)\to Sp(n+1)$. 
Let $Sp:=\colim_{n\to \infty} Sp(n)$.  Recall that 
we have an isomorphism 
$H^{*}(B_{com}Sp(n))\cong M^{B_{n}}/L_{n}$
where $L_{n}=J_{Sp(n)}$ is the ideal in $M^{B_{n}}$  generated 
by the power sums $p_{n}(2,0),\dots, p_{n}(2n,0)$. Thus 
for $Sp$ we have  
\[
H^{*}(B_{com}Sp)\cong \varprojlim H^{*}(B_{com}Sp(n)) \cong 
\varprojlim M^{B_{n}}/L_{n}.  
\]
Define
\[
M^{B_{\infty}}= \varprojlim M^{B_{n}} \text{ and }
L_{\infty}:= \varprojlim L_{n}.
\]
Thus for  the group $Sp$ we have an isomorphism 
of algebras over $\Q$
\[
H^{*}(B_{com}Sp)\cong M^{B_{\infty}}/L_{\infty}.
\]
Next we show that $M^{B_{\infty}}$ is a polynomial algebra.
To see this we first show that the  power sums 
$p_{n}(a,b)=x_{1}^{a}y_{1}^{b}+\cdots+x_{n}^{a}y_{n}^{b}$,
where $a,b$ runs though all non-negative integers such that 
$0<a+b$ and $a+b$ is even generate $M^{B_{n}}$ as a $\Q$-algebra.  
For this consider the averaging operator 
$\rho:\Q[\x,\y]\to \Q[\x,\y]^{^{B_{n}}}$ corresponding to the group 
$B_{n}$. Note that as a $\Q$-vector space 
$M^{B_{n}}$ is generated by 
the elements of the form $\rho(m(\x,\y))$, where 
$m(\x,\y)=x_{1}^{i_{1}}y_{1}^{j_{1}}\cdots x_{n}^{i_{n}}y_{n}^{j_{n}}$ 
is a monomial. By  \cite[Lemma 3.2]{Gomez} if $i_{k}+j_{k}$ is 
odd for some $1\le k\le n$ then $\rho(m(\x,\y))=0$. It follows that 
as a $\Q$-vector space  $M^{B_{n}}$ is generated by 
the elements of the form $\rho(m(\x,\y))$, where 
$m(\x,\y)=x_{1}^{i_{1}}y_{1}^{j_{1}}\cdots x_{n}^{i_{n}}y_{n}^{j_{n}}$  
and $i_{k}+j_{k}$ is even for $1\le k\le n$. Suppose that 
$m(\x,\y)$ is such a monomial. Define the length of $m(\x,\y)$, 
$\ell(m(\x,\y))$, to be the number 
of tuples $(i_{k},j_{k})$ that are nonzero for $1\le k\le n$. 
We can show that 
$\rho(m(\x,\y))$ is a polynomial on the different  $p_{n}(a,b)$ in an 
inductive way on the length of the monomial $m(\x,\y)$.
If 
$\ell(m(\x,\y))=1$ we have 
\[
\rho(m(\x,\y))=\frac{1}{n}\left(x_{1}^{i}y_{1}^{j}+
\cdots x_{n}^{i}y_{n}^{j}\right)
=\frac{p_{n}(i,j)}{n}
\]
and there is nothing to prove. Given any monomial 
$m(\x,\y)=x_{1}^{i_{1}}y_{1}^{j_{1}}\cdots x_{n}^{i_{n}}y_{n}^{j_{n}}$ 
with length $\ell(m(\x,\y))=r$, let 
$(i_{k_{1}},j_{k_{1}}),\dots, (i_{k_{r}},j_{k_{r}})$ be the corresponding 
different tuples that are nonzero. Note that 
\[
p_{n}(i_{k_{1}},j_{k_{1}})\cdots p_{n}(i_{k_{r}},j_{k_{r}})
=\left( \sum_{i=1}^{n}x_{i}^{i_{k_{1}}}y_{i}^{j_{k_{1}}}  \right)\cdots 
\left( \sum_{i=1}^{n}x_{i}^{i_{k_{r}}}y_{i}^{j_{k_{r}}}  \right)
=c\rho(m(\x,\y))+q(\x,\y),
\]
where $c$ is a non-zero constant and $q(\x,\y)$ is a sum 
of certain monomials $n(\x,\y)$ with $\ell(n(\x,\y))<r$. This proves that 
the elements $p_{n}(a,b)$, where $a+b>0$ is even generate 
$M^{B_{n}}$. In fact it can be seen that the collection $\{p_{n}(a,b)\}$, 
where $a,b$ run through all the non-negative integers 
such that $0<a+b\le 2n$ and $a+b$ is even generate $M^{B_{n}}$ but 
we don't need that fact. Recall that  for every $n\ge 1$ we have a map 
$j^{*}_{n}:M^{B_{n+1}}\to M^{B_{n}}$  obtained by sending 
$x_{i}\mapsto x_{i}$, $y_{i}\mapsto y_{i}$  for $1\le i\le n$ and 
$x_{n+1}, y_{n+1}\mapsto 0$ and  
$M^{B_{\infty}}:= \varprojlim M^{B_{n}}$ 
Suppose that $a,b$ are non-negative integers. 
Note that $j^{*}_{n}(p_{n+1}(a,b))=p_{n}(a,b)$ and thus the different 
polynomials $p_{n}(a,b)$ induce an element $p(a,b)$ in 
$M^{B_{\infty}}$. Note that each signed-multisymmetric polynomial 
is in particular a multisymmetric polynomial; that is,  we can 
see $M^{B_{\infty}}$ as a subset of $M^{\Sigma_{\infty}}$. Also, 
we know that $M^{\Sigma_{\infty}}$ is a polynomial algebra on the 
different elements $p(a,b)$  where $a+b>0$ by 
\cite[Theorem 2]{Vaccarino}. This implies in particular that the 
collection $\{p(a,b)\}_{a+b>0, \text{even}}$ is algebraically 
independent in $M^{\Sigma_{\infty}}$ and in particular, it is also 
algebraically independent in $M^{B_{\infty}}$. 
As a corollary the following is obtained.

\begin{corollary}\label{cor signed multisymmetric}
The $\Q$-algebra $M^{B_{\infty}}$ is a polynomial algebra on the 
generators $p(a,b)$, where $a,b$ run through all the non-negative 
integers such that $0<a+b$ is even.  
\end{corollary} 

Using the previous fact we can obtain a description of 
$H^{*}(B_{com}Sp)$ as an algebra. Indeed, recall that 
\[
H^{*}(B_{com}Sp)\cong M^{B_{\infty}}/L_{\infty},
\]
where $L_{\infty}$ is the ideal generated by the power sums 
$p(2n,0)$ for all $n\ge 1$.  Similarly define 
$Sp_{\infty}(\C)=\colim_{n\to \infty} Sp_{n}(\C)$.
As a corollary the following is obtained.

\begin{corollary}\label{comb 3}
Define  
$\Ma_{Sp}=\{(a,b)\in \N^{2}~|~ b>0,\  a+b \text{ even}\}$. 
Then we have isomorphisms of $\Q$-algebras
\[
H^{*}(B_{com}Sp_{\infty}(\C))\cong 
H^{*}(B_{com}Sp)\cong \Q[z_{a,b}~|~ (a,b)\in\Ma_{Sp}].
\]
\end{corollary}
 
 \section{Appendix}

The goal of this appendix is to show that for any Lie group 
$G$ the simplicial space $[B_{com}G]_{*}$ is a proper simplicial 
space. This fact was proved in \cite[Theorem 8.3]{AC} for Lie groups 
$G$ that are closed subgroups of $GL_{n}(\C)$ for some $n\ge 0$ 
and extended in the equivariant setting in \cite{ACG} for compact 
Lie groups. Here we show that the arguments in \cite{AC} 
can be used to proved this result for any Lie group $G$.
\medskip

We start by recalling some basic definitions. 
Recall that a pair  of topological spaces $(X,A)$
is said to be an NDR pair  if 
there exist continuous functions 
$h:X\times  [0,1]\to X$ and $u:X\to [0,1]$
such that the following conditions are satisfied:
\begin{enumerate}
\item $A=u^{-1}(0)$,
\item $h(x,0)=x$ for all $x \in X$,
\item $h(a,t)=a$ for all $a\in A$ and all $t\in [0,1]$, and
\item $h(x,1)\in A$ for all $x\in u^{-1}(\left[0,1\right)).$
\end{enumerate} 
In this case we say that $(h,u)$  is a representation 
of $(X,A)$ as an NDR pair. If in addition we have $u(h(x,t))<1$ 
for all $t\in [0,1]$ 
whenever $u(x)<1$, then $(X,A)$ is called a strong NDR pair. A 
simplicial space  $X_{*}$ is said to be proper if each pair 
$(X_{n+1},sX_{n})$ is a strong NDR pair, 
where $sX_{n}$ is the image of the different 
degeneracy maps in $X_{n+1}$.

\begin{proposition}\label{proper simplicial}
If $G$ is a Lie group then $[B_{com}G]_*$ 
is  a proper simplicial space. 
\end{proposition}

\Proof
Suppose that $G$ is a Lie group 
and let $\g$ denote its Lie 
algebra endowed with a norm $\left\|\cdot\right\|_{\g}$. 
Recall that the exponential map $\exp:\g\to G$ is a 
local homeomorphism. Let $U\subset \g$ be any $Ad$ invariant 
open neighborhood of $0\in \g$ on which the exponential map is 
injective. Fix some  $\epsilon >0$ 
such $\bar{B}_{\epsilon}(0)\subset U$. Then in particular
$\exp_{|}:\bar{B}_{\epsilon}(0)\to G$
is a homeomorphism onto its image. Define a  function
$u:G\to [0,1]$
as follows:
\begin{equation*}
u(g)=\left\{ 
\begin{array}{ll}
2\left\|y\right\|_{\g}/\epsilon & \text{if }g
=\exp (y)\text{ for }g\in \exp (\bar{B}_{\epsilon /2}(0)), \\ 
1 & \text{if }g\in G-\exp (B_{\epsilon /2}(0)).%
\end{array}%
\right. 
\end{equation*}
Also, let  $s:G\to [0,1] $ be any bump function satisfying 
the following conditions
\begin{equation*}
s(g)=\left\{ 
\begin{array}{ll}
1 & \text{if }g=\exp (y)\text{ for }g\in 
\exp (\bar{B}_{\epsilon /2}(0)), \\ 
0 & \text{if }g\in G-\exp (B_{\epsilon}(0)).
\end{array}%
\right. 
\end{equation*}
Finally, define a homotopy 
\[
h:G\times [0,1]\to G
\] 
by 
\begin{equation*}
h(g,t)=\left\{ 
\begin{array}{ll}
\exp((1-t)y) & \text{if }g=\exp (y)\text{ for }
y\in \bar{B}_{\epsilon/2}(0), \\ 
\exp((1-s(g)t)y) & \text{if }g=\exp (y)\text{ for }y\in 
\bar{B}_{\epsilon}(0)-B_{\epsilon /2}(0), \text{ and}\\ 
g & \text{if }g\in G-\exp (B_{\epsilon}(0)).%
\end{array}%
\right. 
\end{equation*}
The functions $h$ and $u$ are defined so that 
$(h,u)$ is representation of $(G,\{1_{G}\})$ as an NDR pair. 
This can be seen in the same way as in \cite[Proposition 8.2]{AC}.
Moreover, we claim that the function $h$ satisfies the following 
property: 
for each $g\ne 1_{G}$ in $G$ and each $0\le t<1$ we have 
$Z_{G}(h(g,t))=Z_{G}(g)$. Indeed, assume that 
$g\in G$ with $g\ne 1_{G}$. Note that if $0\le t<1$ then 
$h(g,t)=g$ if $g\in G-\exp (B_{\epsilon}(0))$ and  
$h(g,t)=\exp(ky)$ for some $0<k\le 1$  if $g=\exp(y)$ with 
$y\in \exp (B_{\epsilon}(0))$. In the first case we have nothing 
to prove. Suppose then that $g=\exp(y)$ with 
$y\in \exp (B_{\epsilon}(0))$ and thus $h(g,t)=\exp(ky)$ for some $0<k\le 1$. 
Since $y\in \exp (B_{\epsilon}(0))\subset U$  and $U$ is an 
$Ad$ invariant open set on which the exponential map is 
injective, then by  \cite[Lemma 3.2.1]{Duistermaat} we have
$Z_{G}(g)=Z_{G}(\exp(y))=Z_{G}(\exp(ky))=Z_{G}(h(g,t))$ 
proving that $Z_{G}(g)=Z_{G}(h(g,t))$ as claimed. 
By \cite[Theorem 7.3]{AC} we conclude that 
the inclusion map 
$I_{j}:S_{n}(j,G)\hookrightarrow S_{n}(j-1,G)$ is a cofibration 
for every $1\le j\le n$.\footnote{In the terminology of \cite[Definition 6.1]{AC}, $G$ has
cofibrantly commuting elements.} Here 
$S_{n}(j,G)$ denotes the subspace of $\Hom(\Z^{n},G)$ consisting 
of the commuting $n$-tuples with at least $j$ coordinates 
equal to $1_{G}$. This implies in particular that 
the inclusion map $s([B_{com}G]_{n-1})=S_{n}(1,G)\hookrightarrow 
\Hom(\Z^{n},G)=[B_{com}G]_{n}$ is a  cofibration.  
Using the explicit NDR-pair representation of 
$(\Hom(\Z^{n},G),S_{n}(1,G))$ provided by \cite[Theorem 7.3]{AC}  
it can easily be seen that $(\Hom(\Z^{n},G),S_{n}(1,G))$ is actually 
a strong NDR-pair, proving that 
$[B_{com}G]_{*}$ 
is  a proper simplicial space. 
\qed
 
\begin{remark}\label{strictly proper}
Using the explicit NDR representation provided 
by \cite[Theorem 7.3]{AC} for the pair  $(\Hom(\Z^{n},G),S_{n}(1,G))$, 
it can be seen that 
in fact  $[B_{com}G]_*$ 
is a  strictly proper simplicial space (see  
\cite[Definition 11.2]{May} for the definition of a strictly proper 
simplicial space). 
\end{remark}

\end{document}